\documentclass{article}
\usepackage{amssymb}
\usepackage{graphicx}
\usepackage{amsmath}

\newtheorem{theorem}{Theorem}[section]

\newtheorem{conjecture}[theorem]{Conjecture}
\newtheorem{corollary}[theorem]{Corollary}

\newtheorem{definition}[theorem]{Definition}

\newtheorem{lemma}[theorem]{Lemma}

\newtheorem{problem}[theorem]{Problem}

\newtheorem{proposition}[theorem]{Proposition}
\newtheorem{remark}[theorem]{Remark}

\newenvironment{proof}{\noindent{\bf Proof.}}{$\square$ \newline \medskip}

\newcommand{\SL}{\mathrm{SL}}
\newcommand{\GL}{\mathrm{GL}}
\newcommand{\EL}{\mathrm{EL}}

\newcommand{\Aut}{\mathrm{Aut}}
\newcommand{\Out}{\mathrm{Out}}
\newcommand{\Mat}{\mathrm{Mat}}

\newcommand{\mr}[1]{\mathcal{#1}}

\newcommand{\la}{\langle}
\newcommand{\ra}{\rangle}

\newcommand{\R}{\mathbb{R}}

\newcommand{\Z}{\mathbb{Z}}

\newcommand{\F}{\mathbb{F}}
\newcommand{\Id}{\mathrm{Id}}
\newcommand{\KC}{\mr{K}}

\newcommand{\KaC}{Kazhdan constant}
\newcommand{\Sym}{\mathrm{Sym}}
\newcommand{\Alt}{\mathrm{Alt}}
\newcommand{\Ind}{\mathrm{Ind}}

\begin{document}
\title{Symmetric Groups and Expander Graphs}
\author{Martin Kassabov}

\maketitle
{
\renewcommand{\thefootnote}{}
}
\begin{abstract}
We construct explicit generating sets $S_n$ and $\tilde S_n$ of the %bounded size for the
alternating and the symmetric groups, which turn the Cayley graphs
$\mr{C}(\Alt(n), S_n)$ and $\mr{C}(\Sym(n), \tilde S_n)$ into a
family of bounded degree expanders for all $n$. This answers
affirmatively an old question which has been asked many times in the
literature. These expanders have many applications in the theory of
random walks on groups, card shuffling
and other areas. %of mathematics.

%{\bf rewrite}
\end{abstract}

\maketitle

{
\renewcommand{\thetheorem}{\arabic{theorem}}
\section*{Introduction}

A finite graph is called an expander if for any (not too big)
set of vertices there are many edges leaving this set. This implies that
expander graphs are highly connected and have
a small diameter. Such graphs have many
practical applications, for example in construction of
computer networks.

Using simple counting arguments it can be shown
that the random $k$-regular graphs are expanders for $k \geq 5$.
However these expanders are not
sufficient for many applications where one needs explicit
families of expander graphs. Constructing such %explicit
examples is a
difficult problem.

A natural candidate for a family of expanders are
the Cayley graphs $\mr{C}(G_i,S_i)$ of a sequence of
finite groups $G_i$ with respect to some (suitably chosen) generating
sets $S_i$. It is known that if there is a uniform bound for the
size of the generating sets $S_i$ then the expanding
properties of the Cayley graphs are related to the representation
theory of the groups $G_i$, more specifically to their
Kazhdan constants.

Using this connection G.~Margulis in~\cite{Mar} gave
the first explicit construction of a family of expanders,
using the Kazhdan property \emph{T} of $\SL_3(\Z)$.
Currently there are several different constructions of expanders
using the representation theory of infinite groups
---
typically one finds a finitely generated infinite group $G$
with a `nice' representation theory (usually the group has
some variant of property \emph{T}, property $\tau$, Selberg
property etc.). In this case the Cayley graphs
of (some) finite quotients $G_i$ of $G$ with respect to the images of
a generating set $S$ of the big group form an expander family.
It is very interesting to ask when can one do the opposite ---
which  leads to the following difficult problem (see~\cite{lu}):

\begin{problem}
\label{makexp}
Let $G_i$ be an infinite family of finite groups. Is it possible
to make their Cayley graphs expanders using suitably chosen
generating sets?
\end{problem}

Currently there is no theory which can give
a satisfactory answer to this question. The answer is known
only in a few special cases: If the family of finite groups
comes from a finitely generated infinite group with
property \emph{T} (or its weaker versions)
then the answer is YES.%
\footnote{The opposite is also true:
For any infinite family of finite groups $G_i$ the existence of generating
sets $S_i$ such that the Cayley graphs $\mr{C}(G_i,S_i)$ are a family of expanders
is equivalent to the existence of a finitely generated subgroup of $\prod G_i$ which has
a variant of property \emph{T} (more precisely property $\tau$ with respect to the
induced topology from the product topology on $\prod G_i$).}
Also if all groups in the family are
``almost'' abelian then the answer is NO (see~\cite{lubwiess})
and this is essentially the only case where a negative answer
of Problem~\ref{makexp}
is known.

A natural family of groups which are sufficiently far from abelian
is the family of all symmetric groups.
The special case of Problem~\ref{makexp} for the symmetric groups, i.e.,
the existence of a generating sets which make thier
Cayley graphs expanders, is
an old open question which has been asked
several times in the literature, see~\cite{BHKLS,expanderbook,lu,LubZuk}.

The asymptotic as $n\to \infty$ of \KaC\ of the symmetric group $\Sym(n)$
with respect to
some natural generating sets are known, see~\cite{BaHarpe}.
Unfortunately in all known examples the \KaC\ goes to zero as the size of
the symmetric groups increase (even though in many cases the sizes
of the generating sets are not bounded), which suggest that
Problem~\ref{makexp} has a negative answer for the family of all symmetric
groups.

On the other hand the symmetric group $\Sym(n)$ can be viewed as
a general linear group over ``the field'' with one element, see~\cite{J2}.
In~\cite{KSL3k}, it is shown that the Cayley graphs of
$\SL_n(\F_p)$ for any prime $p$ and infinitely many $n$ can be
made expanders simultaneously by choosing a suitable generating
sets. Using the previous remark this presents a strong supporting evidence that
Problem~\ref{makexp} has a positive answer.% in the case of symmetric groups.

The main result of this paper%
\footnote{This result was announced in~\cite{Ksym}.}
answers affirmatively Problem~\ref{makexp} in the case of
alternating/symmetric groups.

\begin{theorem}
\label{main}
For all $n$ there exists an explicit  generating set $S_n$
(of size at most $L$) of the alternating group $\Alt(n)$,
such that the Cayley graphs $\mr{C}(\Alt(n),S_n)$
form a family of $\epsilon$-expanders.
Here, $L$ and $\epsilon >0$ are some universal constants.
\end{theorem}

The proof uses the equivalence between family of expanders and
groups with uniformly bounded \KaC s.
Using bounded generation and relative \KaC\ of some small
groups, we can obtain lower bounds for the \KaC s of the symmetric groups
$\Sym(n)$ with respect to several different generating sets.
All these estimates relay on Theorem~\ref{altexp6}, whose
proof uses upper bounds for the characters of the symmetric group and
estimates of the mixing time of random walks.

Theorem~\ref{main} has many interesting applications:
%\begin{itemize}
%\item
First, it provides one of the few constructions of an expander family
of Cayley graphs $\mr{C}(G_i,S_i)$ such that the groups $G_i$
are not obtained as quotients of some infinite group having a variant
of Kazhdan property \emph{T}.%
\footnote{As mentioned before if $\mr{C}(G_i,S_i)$ are expanders, then
there exists an infinite group with a variant of property \emph{T}.
 The main point here is that
we prove that the Cayley graphs are expanders without using the representation
of this infinite group.}
The other constructions which do not use a variant of property
\emph{T} are based on an entirely different idea ---
the so called `zig-zag' product of graphs,
for details see~\cite{ALW,MW,RVW,RSW}.

%\item
Second, the automorphism groups of the free
group $\Aut(\mr{F}_n)$ can
be mapped onto infinitely many alternating groups, see~\cite{Gilman}. This, together with
Theorem~\ref{main} provides a strong supporting evidence
for the conjecture that $\Aut(\mr{F}_n)$ and $\Out(\mr{F}_n)$
have property $\tau$.
This conjecture if correct, will imply that the product replacement algorithm
has a logarithmic mixing time, see~\cite{LP} for details.

%\item
Third, Theorem~\ref{main} implies that for a fixed $C$,
the expanding constant of
$\Alt(n)$ with respect to the set $S_n^{C}$ is large enough.%
\footnote{More precisely we have that the spectral gap of the normalized Laplacian of
the Cayley graph is very close to $1$.}
The size of the set $S_n^{C}$ is independent on $n$, and
if $n$ is sufficiently large then $|S_n^{C}| < 10^{-30}n^{1/30}$.
The last inequality allows us to use the expander
$\mr{C}(\Alt(n);S_n^{C})$ as a
`seed' graph for
recursive construction of expanders suggested by
E.~Rozenman, A.~Shalev and A.~Widgerson
in~\cite{RSW}.
This will be one of the few constructions of an infinite family of
expander graphs which are purely combinatorial, i.e., it does
not use any representation theory.
This construction produces a family of expander graphs
from the automorphism groups of
$n$-regular rooted tree of depth $k$. A slight modification of
this construction gives another recursive expander family based on
$\Alt(n^k)$ for fixed large $n$ and different $k$-s.

%\end{itemize}

Theorem~\ref{main} implies the analogous result for the symmetric groups:
\begin{theorem}
\label{mainsym}
For all $n$ there exists an explicit  generating set $\tilde S_n$
(of size at most $L$) of the alternating group $\Sym(n)$,
such that the Cayley graphs $\mr{C}(\Sym(n),\tilde S_n)$
form a family of $\epsilon$-expanders.
Here $L$ and $\epsilon >0$ are some universal constants.
\end{theorem}

%\begin{itemize}
%\item
Theorem~\ref{mainsym} implies that the random walk on
$\Sym(n)$,
generated by $\tilde S_n$
% respectively,
has mixing time approximately
%\linebreak
$O(\log |\Sym(n) |) =O( n \log n)$ steps.
\iffalse
This leads to a natural and
fast algorithm for generating pseudo-random permutations.
There are others algorithms which generate a random permutations in
$\Sym(n)$ in $n \log n$ steps; however they usually require generating
a random number from $1$ to $n$ on each step.
\fi

%\item
%{\bf Fast card shuffles}
%\end{itemize}
%\newpage

%\noindent
The rest of the paper is organized as follows:
Section~\ref{outline} contains definitions and
a sketch of the proof of Theorem~\ref{altexp6}, which is a
weaker version of Theorem~\ref{main}.
The detailed proof of this theorem is contained in
Sections~\ref{relT} and~\ref{cube}.
Section~\ref{proofmain} explains how Theorems~\ref{main} and~\ref{mainsym}
can be derived from Theorem~\ref{altexp6}.
Section~\ref{comments} concludes the paper with some
comments about possible modifications and applications
of Theorem~\ref{main}.

\medskip

\noindent
\textbf{Acknowledgements: }
I wish to thank Alex Lubotzky and Nikolay Nikolov %, Y.~Roichman, E.~Rozemman
for their encouragement and  useful discussions during the work on this project.
I am also very grateful to Yehuda Shalom and Efim Zelmanov for introducing me to the subject.

}

\section{Outline}
\label{outline}

Let us start with one of the equivalent definitions of expander graphs:

\begin{definition}
\label{expander}
A finite graph $\Gamma$ is called an $\epsilon$-expander for some
$\epsilon \in (0,1)$ if for any subset $A \subseteq \Gamma$ of size
at most $|\Gamma|/2$ we have $|\partial (A)| >\epsilon |A|$.
Here $\partial(A)$ is the set of vertices of $\Gamma \backslash A$ of edge
distance 1 to $A$. The largest such $\epsilon$ is called the
expanding constant of $\Gamma$.
\end{definition}

Constructing families of $\epsilon$-expanders
with a large expanding constant $\epsilon$ and bounded valency
is an important practical problem in computer science,
because the expanders can be used to construct concentrators,
super concentrators, contractors and etc.
For an excellent
introduction to the subject we refer the reader to the
book~\cite{expanderbook} by A.~Lubotzky and to~\cite{Kl}.
%{\bf other references?}

If we consider only graphs $\Gamma_i$
were the degree of each vertex is at most $k$, then
all graphs $\Gamma_i$ are $\epsilon$-expanders for some $\epsilon$, if
any of the following equivalent conditions holds:%
\footnote{If the degree of the graphs is not bounded, there are two
different notions of expander graphs -- one corresponding to
Definition~\ref{expander} and a bound of the Chegeer constant; and a second
one, coming from a bound on the spectral gap of the Laplacian.
In the rest of the paper we will use
the second definition, which is more restricitve.}
\begin{enumerate}
\item
The Chegeer constants of $\Gamma_i$ are uniformly bounded away from zero;

\item
The Laplacian of $\Gamma_i$ have a uniformly bounded spectral gap;
\end{enumerate}
In this case we also have that
the lazy random walk on the graph
$\Gamma_i$ mixes in $O(\log |\Gamma_i|)$ steps.

%It is well known that the random graphs of degree at least  $3$
%are expanders, but for many applications one needs explicit graphs.
The graphs which appear in many applications arise from finite groups ---
they are Cayley with respect to some generating set or their quotients.
In that case we have an additional equivalent condition --- the
Cayley graphs $\mr{C}(G_i;S_i)$ of $G_i$,
with respect to generating sets $S_i$ % of bounded size,
are $\epsilon$-expanders for some $\epsilon$ if and only if
the Kazhdan constants $\mr{K}(G_i;S_i)$ are uniformly bounded away from zero.%
%\footnote{
%If one considers some  quotients of the Cayley graphs (known as Schreier graphs)
%then there is an analogues condition involving the relative \KaC.}
%Here we use the following definition of property \empf{T} and the Kazhdan Constant:

\medskip

The original definition of Kazhdan property \emph{T} uses the Fell
topology of the unitary dual of a group, see~\cite{kazhdan}.
We are interested not only in property \emph{T}, which
automatically holds for any finite group,
%trivial for finite groups,
but also in the related notion of \KaC s.
The following definitions, which addresses the
notion of the \KaC s,
are equivalent %(only in the case of discrete groups)
to the usual definitions of relative property \emph{T} and
property \emph{T}.

\begin{definition}
\label{tdef}
Let $G$ be a discrete group generated by a finite set $S$ and
let $H$ be a subset of $G$.
Then the pair $(G,H)$ has relative property \emph{T}
if there exists $\epsilon  >0$
such that for every unitary
representation $\rho: \ G \rightarrow U(\cal{H})$ on a
Hilbert space $\mr{H}$ without $H$ invariant vectors and
every  vector $v \not = 0$ there
is some $s \in S$ such that $||\rho(s)v-v||> \epsilon ||v||$.
The largest $\epsilon$ with this property is called the
\emph{relative \KaC} for $(G,H)$ with respect to the set
$S$ and is denoted by $\KC(G,H;S)$.

The group $G$ has Kazhdan property \emph{T} if the pair $(G,G)$ has
relative \emph{T} and the \emph{\KaC} for the group $G$ is
$\KC(G;S) := \KC(G,G;S)$.
\end{definition}

The property \emph{T} depends only on the group $G$
and does not depend on the choice of the
generating set $S$, however the Kazhdan constant
depends also on the generating set.

It is clear that any finite group $G$ has property \emph{T}, because
it has only finitely representations generated by a single vector. % ---
%without loss of generality in Definition~\ref{tdef} we
%can requite that $\mr{H}$ is generated by the vector $v$.
If the generating set $S$ contains all elements of the group $G$,
then we have the inequalities
$$
2 \geq \KC(G;G)\geq \sqrt{2}.
$$
This follows from the following observation:
if a unit vector $v \in \cal{H}$
is moved by at most $\sqrt{2}$ by any element of the group $G$,
then the whole orbit $Gv$ is contained in some half space and its
center of mass is not zero. The $G$ invariance of the orbit
gives that the center of mass  is a non-zero $G$-invariant
vector in $\cal{H}$. However the resulting Cayley graphs are complete graphs
with $|G|$ vertices and very `expensive' expanders.

\iffalse
\medskip

Some estimates for \KaC\ of the symmetric groups with respect to
some generating sets are known ---
from~\cite{Ro1} it follows that
$$
\KC(\Sym(n);T_{k,n}) \preceq \frac{k}{\sqrt{n}},
$$
where $T_{k,n}$ is the set of all permutations with less that
$k$ non-fixed points. This implies that
the Cayley graphs $\mr{C}(\Sym(n);T_{k,n})$
are not expanders if $k \ll \sqrt{n}$.
The size of the generating set $T(k,n)$ in this
case is $T(k,n) \approx n^{\sqrt{n}}$. We must mention
that so far explicit expanding
generating sets of the symmetric group which are significantly
smaller then the symmetric group itself, are not known ---
the size of the smallest
one is around $n^{\sqrt{n}}$, in~\cite{Ro1} it is shown that
$\mr{C}(\Sym(n);T_{k,n})$ have some expanding properties
if $k \gg \sqrt{n}$.

On the other hand, it is known
%(Alon-Roichman Theorem~\cite{AR},
%and its slight improvements~\cite{LR} and~\cite{LS})
that
the Cayley graphs $\mr{C}(\Sym(n);R)$ are expanders if $R$ is a
random generating set of size $n\log n$. This follows form
a general result, known as Alon-Roichman Theorem~\cite{AR}
%(and its slight improvements~\cite{LR} and~\cite{LS}),
(see also~\cite{LR} and~\cite{LS}),
valid for all finite groups. This result is the best possible
in the class of all finite groups, however it is believed
that in the case of the symmetric group the bound $n\log n$
can be improved significantly.

\medskip
\fi

As mention before to prove Theorem~\ref{main}, it is enough to
prove that there exist generating sets, the \KaC s of which
are uniformly bounded away from $0$.
We will start with a similar result for $\SL_n(\F_p)$ for a fixed $n\geq 3$,
which also proves that $\SL_n(\Z)$ has property \emph{T}
and gives an estimate for the \KaC\  with respect to the generating
set consisting of the elementary matrices.

The standard proofs of property \emph{T} for arithmetic groups, like
$\SL_n(\Z)$, use the representation theory of Lie groups and are not
quantitative, i.e., they do not lead to any estimate of the \KaC s.
Our approach is based on ideas from~\cite{YSh}: we start with
the trivial estimate
$$
\KC(\SL_n(\F_p);\SL_n(\F_p)) \geq \sqrt{2}
$$
and  will make several
changes of the initial generating set $\SL_n(\F_p)$, which will decrease
its size and keep some estimate of the \KaC. Finally we will
end with a generating $S_{n,p}$ of size independent on $p$ such that
$\KC(\SL_n(\F_p);S_{n,p})$ is bounded away from $0$.

There are two propositions which allow us to estimate the
\KaC \
if we change the generating set. The first one is a quantitative
version of the fact that
property \emph{T} for a discrete group does not depend on the generating set.

\begin{proposition}
\label{KCball}
Let  $S$ and $S'$ be two finite generating sets of a group $G$,
such that $S' \subset S^k$, i.e., all the elements of $S'$ can
be written as short products of elements form $S$. Then we have
$$
\KC(G;S) \geq \frac{1}{k}\KC(G;S').
$$
\end{proposition}
\begin{proof}
If $v$ is an $\epsilon$ almost invariant vector for the set $S$ then
$$
{\renewcommand{\arraystretch}{2}
\begin{array}{r@{}c@{}l}
||\rho(s_1\dots s_k) v -v || & {}= {} & %\displaystyle
\left|\left|
\sum_{j=1}^k \rho(s_1\dots s_{j-1})\left(\rho(s_j) v -v\right)
\right|\right| \leq \\
& \leq & %\displaystyle
\sum_{j=1}^k ||\rho(s_j) v -v|| \leq k\epsilon.
\end{array}
}
$$
This shows that  $v$ is $k\epsilon$ almost invariant
vector for the set $S^k$ and in particular for $S'$.
If
we start with $\epsilon = \frac{1}{k}\KC(G;S')$
then $v$ is $\KC(G;S')$
almost invariant vector for $S'$, which gives the
existence of nonzero  invariant vectors in $\mr{H}$.
\end{proof}

\medskip

The second proposition uses relative property \emph{T} to enlarge the
generating set by adding a subgroup to it:
\begin{proposition}
\label{KCrel}
If $H$ is a normal subgroup of a group $G$ generated by a set $S$ then:
$$
\KC(G;S) \geq \frac{1}{2}\KC(G,H;S)\KC(G;S \cup H).
$$
\end{proposition}
\begin{proof}
Let $\rho: G \to U(\mr{H})$ is a unitary representation and let
$v$ be $\epsilon$ almost invariant vector for the set $S$.
We can write
$\mr{H} = \mr{H}^{||} \oplus \mr{H}^{\perp}$,
where $\mr{H}^{||}$ is the space of all $H$ invariant vectors in $\mr{H}$
and $\mr{H}^{\perp}$ is its orthogonal complement.
These spaces are
$G$-invariant because $H$ is a normal subgroup.
This decomposition gives $v= v^{||} + v^{\perp}$, where both
components are $\epsilon$-almost invariant vectors. However there are
no $H$ invariant vectors in $\mr{H}^\perp$
and relative property \emph{T} of $(G,H)$ implies
%therefore
$$
||\rho(s) v^\perp - v^\perp|| \geq \KC(G,H;S)||v^\perp||
$$
for some $s\in S$, which implies that
$
||v^\perp|| \leq \epsilon \KC(G,H;S)^{-1}.
$
Thus, for any $h\in H$ we have
$$
||\rho(h) v - v || = || \rho(h) v^\perp - v^\perp|| \leq 2 ||v^\perp|| \leq
2\epsilon \KC(G,H;S)^{-1}.
$$
If we start with $\epsilon = \frac{1}{2}\KC(G,H;S)\KC(G;S \cup H)$ than
the above inequality gives that $v$ is $\KC(G;S \cup H)$-almost
invariant for both $H$ and $S$
therefore there exists an invariant vector in $\mr{H}$.
\end{proof}

In practice one often uses a variant of this proposition for several
subgroup $H_i$ simultaneously.
\begin{proposition}
\label{relKC}
a) Let $H_i$ and $N_i$ be subgroups of a group $G$ such that
\linebreak
$H_i \lhd N_i < G$ and $N_i$ is generated by $S_i$.
If $\KC(N_i,H_i;S_i) \geq \alpha$
and $S_i \subset S$ for all $i$  then
$$
\KC(G;S) \geq \frac{1}{2} \alpha \KC(G; \cup H_i).
$$

b) Let $N$ be a group generated by a set $\overline{S}$ and let $H \lhd N$.
If $\pi_i$ are homomorphisms from $N$
in a group $G$, then
$$
\KC(G;\cup \pi_i(\overline{S})) \geq \frac{1}{2} \KC(N,H;\overline{S}) \KC(G; \cup \pi_i(H)).
$$
\end{proposition}

\medskip

Using these simple propositions Y.~Shalom~\cite{YSh} was able to estimate the
Kazhdan constant of $\SL_n(\Z)$ with respect to the set $S$
of all elementary matrices with $\pm 1$ off the diagonal.
Let $EM_{i,j}$ denote the subgroup of elementary matrices of the
form $\Id + n e_{i,j}$ and let $EM = \bigcup_{i\not=j} EM_{i,j}$.

The relative \KaC\ of the pair $(\SL_2(\Z) \ltimes \Z^2,\Z^2)$ was
estimated by Burger~\cite{Bur}, see also~\cite{GG2}, to be:
%(Burger attributes the idea of the bound to Furstenburg??)
$$
\KC(\SL_2(\Z) \ltimes \Z^2,\Z^2;\overline{S}) \geq 1/10,
$$
where $\overline{S}$ is the set consisting of the
$4$ elementary matrices in $\SL_2(\Z)$ together with the
standard basis vectors of $\Z^2$ and their inverses.
The proof is a quantitative version of the fact that there
are only a few $\SL_2(\Z)$ invariant measures on the torus $\R^2/\Z^2$.

Using embeddings $\pi_{i,j}$ of
$\SL_2(\Z) \ltimes \Z^2$ in $\SL_n(\Z)$ such that
the image of $\Z^2$ under $\pi_{i,j}$ contains $E_{i,j}$, by
Proposition~\ref{relKC} we have
$$
{\renewcommand{\arraystretch}{2}
\begin{array}{r@{}c@{}l}
\displaystyle
\KC(\SL_n(\Z),S)  & {}\geq {} &
\displaystyle
 \frac{1}{2} \KC(\SL_2(\Z) \ltimes \Z^2,\Z^2;\overline{S})
\KC(\SL_n(\Z);\cup EM_{i,j}) \geq
\\
& \geq &
\displaystyle
\frac{1}{20}\KC(\SL_n(\Z),EM).
\end{array}
}
$$

The next step in the proof uses bounded generation
of $\SL_n(\Z)$
with respect to the set of the elementary matrices.%
\footnote{The bounded generation of $\SL_n(\Z)$ is a deep
result due to Carter and Keller, see~\cite{CK} and~\cite{CK1}.
For our purposes we need only the
\KaC\ of the finite groups $\SL_n(\Z/s\Z)$ and we need
bounded generation of these groups, which follows
immediately from the standard row reduction algorithm.
Similar result is also valid for the finite quotients of
the groups $\SL_n(\Z[x_1,\dots,x_k])$, although the
bounded generation for these groups is not known, see~\cite{KNtau}.}
%This property is equivalent to:
%$EM^N = \SL_n(\Z)$ for $N \geq \frac{3}{2}n^2+ 60$.
This implies that
$$
\KC(\SL_n(\Z);EM) \geq \frac{1}{N}\KC(\SL_n(\Z);EM^N) =
\frac{1}{N}\KC(\SL_n(\Z);\SL_n(\Z))\geq\frac{\sqrt{2}}{N},
$$
where $N \geq \frac{3}{2}n^2+ 60$.

Combining these two inequalities one obtains a lower bound
for the \KaC\
$$
\KC(\SL_n(\Z),S) > \frac{1}{30n^2+1200}.
$$
In particular, this implies that
there is a bound for the \KaC s $\KC(\SL_n(\F_p);S)$
which is independent on $p$, implying that
the Cayley graphs $\mr{C}(\SL_n(\F_p);S)$ of the
finite groups $\SL_n(\F_p)$ are expanders.
\footnote{
In~\cite{K}, it is shown that $\left(\KC(\SL_n(\Z),S)\right)^{-1} = O(\sqrt{n})$.
This yields an asymptotically exact estimate for the expanding
constant of the Cayley graphs of
$\SL_n(\F_p)$ with respect to the set of all elementary
matrices.}
%
%It is interesting to note that if the size of the matrices
%increases, then the resulting Cayley graphs do not form an
%expander family, even though the degree of these graphs
%goes to infinity.

Using relative property \emph{T} of the pair
$\SL_2(R)\ltimes R^2, R^2$ for finitely generated noncommutative rings
$R$, the Cayley graphs of $\SL_n(\F_q)$ for any prime power $q$
and infinitely many $n$ can be
made expanders simultaneously by choosing a suitable generating
set, see~\cite{KSL3k}. An important building block
in this construction is that the group $\SL_n(\F_q)$
can be written as a product of $20$ abelian subgroups and this number is
independent on $n$ and $q$.

\medskip

These methods can not be applied to
the symmetric/alternating groups. Most of the
estimates for the relative \KaC\ $\KC(G,H;S)$
use that there are no invariant measures on
the dual $\widehat H$
of the group $H$ under the action of the normalizer of $H$ in $G$.
The quantitative versions of this fact are known only if
$\widehat H$ is well understood, which is the case
only if the subgroup $H$ is (or is very close to)
an abelian subgroup.

If we start with  finite
generating sets $S_n$ of bounded size of the alternating groups $\Alt(n)$,
then we can find
only finitely many abelian groups $H_\alpha$
(the number depending only on the size of the generating set $S_n$)
such that the relative \KaC s can be estimated easily.
Thus would allow us to bound $\KC(\Alt(n);S_n)$ with
$\KC(\Alt(n);E_n)$, where $E_n$ is a union of a bounded number
of abelian subgroups. We can use Proposition~\ref{KCball}
to estimate $\KC(\Alt(n);E_n)$ if $E_n^k = \Alt(n)$ for some $k$, i.e.,
if each alternating group $\Alt(n)$ is a product of a fixed number
of abelian subgroups.

However, the finite%
\footnote{It is interesting that the full symmetric group on
an infinite set can be written as a product of $250$ abelian
subgroups, see~\cite{A}.}
symmetric/alternating groups do not have this property ---
the size of  $\Sym(n)$ or $\Alt(n)$ is
approximately $n^{n}$ and every abelian subgroup
has no more than $2^n$ elements, thus one needs at least $\ln n$
subgroups.
This suggests that $\Alt(n)$ are ``further from the abelian groups''
than all other finite simple groups, and therefore
they should have more expanding properties.
Unfortunately, this also
significantly complicates the construction of expanders
based on the alternating groups, because the above method
can not be applied %in this case
without significant modifications.%
\iffalse
\footnote{May be if one make several extra steps, it is possible to
obtain a expanding generating sets of the symmetric groups.
We think that this is very unlikely to work because
if one tries to find another system of abelian groups $K_\alpha$
such that $\KC(G_i,K_\alpha;E^k)$ are bounded form bellow then
many of $K_\alpha$ are going to be conjugate by elements in $E^k$,
which will make the size of $\cup K_\alpha$ similar to $E^k$.}
\fi

\medskip

The main idea of the proof of Theorem~\ref{main} is a
modification of the above method --- the difference
is that we are looking for a set $E_n$, which is a union of
finitely many abelian subgroups in $\Sym(n)$, such that the
\KaC s $\KC(\Sym(n);E_n)$ are uniformly bounded, even though there is
no $k$ such that $E_n^k = \Sym(n)$ for all $n$.
There is a natural construction of the set $E_n$ if $n$ has a specific form.
If $n=(2^{3s}-1)^6$ for some $s$, we shall construct a group
$\Delta$ generated by a set $\bar S$
and an abelian subgroup $\Gamma$ in it, such that
the relative \KaC s $\KC(\Delta,\Gamma;\bar S)$ are
bounded away from $0$.
There are several natural
embeddings $\pi_i : \Delta \to \Alt(n)$
and the \KaC\
$\KC(\Alt(n); \cup \pi_i(\Gamma))$
can be estimated.
These two
bounds give us a bound for $\KC(\Alt(n); \cup \pi_i(S))$, which will prove
the following:

\begin{theorem}
\label{altexp6}
If $N=(2^{3s}-1)^6$ for some $s > 6$
there exists a generating set $S_N$
(of size at most $200$) of the alternating group $\Alt(N)$,
such that the Cayley graphs $\mr{C}(\Alt(N),S_N)$
form a family of $\epsilon$-expanders.
Here $\epsilon$ % =2.10^{-5}$
is a universal constant.%
%\footnote{The size of the generating set $S_N$ can be significantly
%decreased to $|S_N|<10$
%at the expense of slightly decresing the expanding constant.
%}
\end{theorem}

\emph{Sketch of the proof of Theorem~\ref{altexp6}}. This sketch
explains the main idea of the proof ---
the complete proof is in Sections~\ref{relT} and~\ref{cube}.
We will think that the alternating group $\Alt(N)$ acts
on a set of $N$ points which are arranged into $d=6$ dimensional
cube of size $K=2^{3s}-1$. The group $\Gamma$
is a direct product of $K^{d-1}=K^5$ cyclic groups of order $K$.
This group can be  embedded in $\Alt(N)$ in $d$ different ways
as follows: the image of each cyclic group
in $\Gamma$ under $\pi_i$ permutes the points on a line, %also celled needle,
parallel to the $i$-th coordinate axis.
An other way to think about $E_i=\pi_i(\Gamma)$ is as part of the
subgroup of $\Alt(N)$ which preserves all coordinates but the $i$-th one.

The group $\Gamma$ can be embedded in a group $\Delta$
generated by a set $\overline{S}$, such that the embeddings $\pi_i$
can be extended to $\Delta$.
The group $\Delta$ is a product of many copies of
$\SL_{3s}(\F_{2}) = \EL_3(\Mat_s(\F_2))$
and can be viewed as $\EL_3(R)$, where the ring $R$ is a product of
many copies of the matrix ring $\Mat_s(\F_{2})$.
An important observation is that the ring $R$ has a generating
set whose size is independent on $s$. Using results
from~\cite{KSL3k}, mainly the relative \KaC\
$$
\KC(\EL_2(\Z \la x_1,\dots,x_k \ra)\ltimes \Z\la x_1,\dots,x_k \ra^2,
\Z\la x_1,\dots,x_k\ra^2; F)
$$
for some set $F$. This allows us to obtain an estimate for $\KC(\Delta;\bar S)$
\begin{equation}
\label{KCGamma}
\KC(\Delta,\Gamma;\overline{S}) \geq
\KC(\Delta;\overline{S})\geq \frac{1}{550}%??
\end{equation}
and to compare $\KC(\Alt(N); \cup \pi_i(\overline{S}))$ and
$\KC(\Alt(N); \cup \pi_i(\Gamma))$.

\medskip

We will finish the proof of Theorem~\ref{altexp6} using
Theorem~\ref{KCGammas} from Section~\ref{cube}, which gives us that
\begin{equation}
\label{KCSym}
\KC\left(\Alt(N); \cup E_i \right) \geq \frac{1}{70},%??,
\end{equation}
where $E_i=\pi_i(\Gamma)$, provided that $s > 6$.%
\footnote{If $s \leq 6$, using the same methods
we can also obtain a bound for the \KaC \ like
$\KC(\Alt(N); \cup \pi_i(\Gamma)) \geq 1/1000$,
but this requires more careful analysis.}
The proof of this inequality uses directly the representation
theory of the alternating group. Let
$\rho$
be a unitary representation of $\Alt(N)$ in a Hilbert space $\mr{H}$
and let $v \in \mr{H}$ be
$\epsilon$-almost invariant vector with respect to
the set $E$. We want to prove that if
$\epsilon < 1/70$ then $\mr{H}$ contains an invariant vector.

\medskip

First, we will split the representation $\rho$ into
two components --- one corresponding
to partitions $\lambda$ with $\lambda_1 < N - h$ and
a second one, containing all other partitions
(here $h$ depends on $N$ and will be determined later).
This decomposition of the representation $\rho$ %of $\Alt(N)$
into two
components $\mr{H}_1$ and $\mr{H}_2$,
using the first part of the partition $\lambda$
borrows ideas from~\cite{Ro1}. There, Y.~Roichman uses a similar
argument to show that the Cayley
graphs of the alternating group with respect to a conjugancy
class with a large number of
non-fixed points have certain expanding properties.

First we shall prove that the projection $v_1$ of the vector
$v$ in the first component $\mr{H}_1$
is small provided that $h \gg K \ln K$.
%({\bf may be even $h<50K$??}),
%however we shell prove a weaker version with $g \gg K^{5/4}$.
It can be shown that for a fixed $k$ the ball
$E^k$ almost contains an
entire conjugancy class $C$ containing all  cycles of
length approximately $K^{d-1}/3\ln K$.
Using estimates for the values of characters of the symmetric group
we can show that if $h \gg K\ln K$ then
$$
\left| \left| \frac{1}{C}\sum_{g\in C} gv_1 \right| \right | \ll ||v_1||.
$$
The vectors $v$ and $v_1$
are almost invariant with respect to the set $\cup \pi_i(\Gamma)$ and
therefore with respect to $C$. This implies that
\begin{equation}
\label{h1o}
\left| \left| \frac{1}{C}\sum_{g\in C} gv_1 - v_1\right|\right|
\leq 47 \epsilon.%??
\end{equation}
The two inequalities gives that $||v_1|| \leq 47 \epsilon + 0.07$.

\medskip

That projection $v_2$ of $v$
in the second component $\mr{H}_2$ is close to an
invariant vector if $h < K^{d/4}$. %
%again we will slightly weaker version only for
%$h \ll K^{d-8/5} =K^{7/5}$
\iffalse
\footnote{We believe that the argument also works even in the case
$h \ll K^{d-3/2}$,
however we are unable to prove it.
If such generalization is true, it will allow us to use $d=3$,
which will improve the estimates
for the \KaC s.% by a factor of $10$.
}
\fi
%
Here the idea is that the space $\mr{H}_2$
can be embedded into a vector space with a basis $\mr{B}$ consisting all
ordered tuples of size $h$. It can be shown that the mixing
time of the random walk on the set $\mr{B}$, generated by
$E$ has
a small mixing time --- this is possible in part because
$|E|\gg |\mr{B}|$, but the proof uses the specific
structure of the set $E$. This implies that the vector
$$
U^K v_2,
\qquad
\mbox{where}
\qquad
U = \frac{1}{|{E}|}\sum_{g\in E} g
$$
is very close to some invariant vector $v_{0}$.
On the other hand this
vector is close to $v_2$, therefore
\begin{equation}
\label{h2o}
||v_2 - v_{0} || \leq 16 \epsilon. %??
\end{equation}

Equations~(\ref{h1o}) and~(\ref{h2o})
were obtained using the  assumptions $h\gg K\ln K$ and $h < K^{d/4}$ ---
these restrictions can be satisfied only if
$K^{d/4} \gg K \ln K$, i.e., $d> 4$.
In order to simplify the argument we require that $d$ is even,
which justifies our choice of $d=6$ and $N=K^6$.
In this case, $h$ have to satisfy the inequalities
$K\ln K \ll h < K^{3/2}$, therefore we chose
$h=\frac{1}{2}K^{3/2}$
and define $\mr{H}_1$ to be the sub-representation of
$\mr{H}$ corresponding to all partitions
with $\lambda_1 < N - \frac{1}{2}K^{3/2}$.

Equations~(\ref{h1}) and~(\ref{h2}) give us that
$$
||v - v_{0}|| \leq 63 \epsilon+ 0.07.%??
$$
In particular this implies that $v_{0} \not =0$ if $\epsilon$ is
small enough and the representation
$\mr{H}$ contains an invariant vector and finishes the proof of
Theorem~\ref{KCGammas}.

Finally equations~(\ref{KCGamma}) and~(\ref{KCSym}) together with
proposition~\ref{relKC} imply that
$$
\KC(\Alt(N);\cup S_i) \geq \frac{1}{2} \KC(\Alt(N); \cup \Gamma_i)
\KC(\Delta,\Gamma;\bar S) > 10^{-5},%??
$$
which completes the proof of Theorem~\ref{altexp6} modulo
the results in Sections~\ref{relT} and~\ref{cube}.
In section~\ref{proofmain},
we derive Theorems~\ref{main} and~\ref{mainsym} from Theorem~\ref{altexp6}.

\section{The groups $\Gamma$ and $\Delta$}
\label{relT}

As mention before, we will think that the
alternating group $\Alt(N)$ acts
on a set of $N$ points which are arranged into $d=6$ dimensional
cube of size $K=2^{3s}-1$ and we will identify these
points with ordered
$6$-tuples of nonzero elements from the field $\F_{2^{3s}}$.
For the rest of this section we will assume that $s$, $K$ and $N$ are
fixed and that $s>6$. For any associative ring $R$, $\EL_3(R)$ denotes the subgroup of the
$3\times 3$ invertible matrices with entries in $R$, generated
by all elementary matrices.

\medskip

Let $H$ denote the group $\SL_{3s}(\F_{2}) = \EL_3(\Mat_s(\F_2))$.
The group $H$
has a natural action on the set $V\setminus \{0\}$ of $K$ nonzero
elements of a vector space $V$ of dimension $3s$ over $\F_{2}$.
The elements of $H$ act by even permutations on
$V\setminus \{0\}$, because
$H$ is a finite  simple group and does not have $\Z/2\Z$ as a factor.

We can identify $V$ with the filed $\F_{2^{3s}}$ --- the
existence of a generator for the multiplicative group of
$\F_{2^{3s}}$ implies
that some element of $H = \GL_{3s}(\F_2)$ acts as a $K$-cycle on
$V\setminus\{0\}$.
We have the equality $\GL_{3s}(\F_2) =\SL_{3s}(\F_2)$, because
we work over a field characteristic $2$. This requirement is
not necessary and can be avoided if we modify the proof of
Theorem~\ref{KCGammas} for other possible choices of
the group $H$, which also give bounded degree expanders in
infinitely many alternating groups see Section~\ref{comments}.

Let $\Delta$ be the direct product of $K^{d-1}$ copies of the group $H$.
The group $\Delta$ can be embedded into $\Alt(N)$ in $6$ different ways
which we denote by $\pi_i$, $i=1,\dots , d$. The image of
each copy of $H$ under $\pi_i$ acts as $\SL_{3s}(\F_2)$
on a set of $K=2^{3s}-1$
points where all
coordinates but the $i$th one are fixed. The existence of
an element of order $K$ in $H$ implies that
$\Delta$ contains an abelian subgroup $\Gamma$
isomorphic to $(\Z/K\Z)^{\times K^{d-1}}$.

\medskip

Another way to think of this group is as follows ---
$\Delta$ is a product of copies of $\SL_{3}(\F_{2^s})$, i.e.,
$$
\Delta \simeq \SL_{3s}(\F_{2})^{\times K^{d-1}} \!\!\!\simeq
\EL_{3}(\Mat_s(\F_{2}))^{\times K^{d-1}} \!\!\! \simeq
\EL_3\left(\Mat_s(\F_2)^{\times K^{d-1}}\right),
$$
i.e., $\Delta \simeq \EL_3(R)$ where $R$ denotes the product of $K^{d-1}$ copies of the
matrix ring $\Mat_s(\F_{2})$.

\begin{lemma}
For any $s$ the ring $R$ is generated by
$2 + \lceil 3(d-1)/s \rceil \leq  5$ elements.
\end{lemma}
\begin{proof}
The matrix algebra $\Mat_s(\F_2)$ can be generated as a ring by
$1$ and two elements $\bar \alpha$ and $\bar \beta$,
for example we can take
$\bar \alpha=e_{2,1}$ and $\bar\beta = \sum e_{i,i+1}$.
By construction, the ring $R$ is
$$
R = \Mat_s(\F_2)^{\times K^{d-1}}
$$
We can think that the copies of $\Mat_s(\F_2)$ are indexed
by tuples of length
$$
t = \lceil \log_{|\Mat_s(\F_2)|} K^{d-1} \rceil  =
\left\lceil \frac{(d-1)\log_2 K }{s^2} \right\rceil \leq
\left\lceil \frac{3(d-1)}{s}\right\rceil
$$
of elements in $\Mat_s(\F_2)$.
Let us define the elements $\alpha$, $\beta$ and $\gamma_i$, $i=1,\dots, t$ in
$\Mat_s(\F_2)^{\times K^5}$ as follows:
All components of $\alpha$ and $\beta$ are equal to
$\bar \alpha$ and $\bar \beta$, the components of $\gamma_i$
in the copy of $\Mat_s(\F_2)$
index by $(p_1,\dots,p_t)$ is equal to $p_i$.
It is no difficult to show that the elements  $\alpha$, $\beta$ and  $\gamma_i$ generate
$\Mat_s(\F_2)^{\times K^{d-1}}$ as associative  ring.
\end{proof}

Any generating set of a ring $R$ gives a generating set of
the group $\EL_3(R)$:
\begin{corollary}
For any  $s$ the group $\Delta$ can be generated a by set $\bar S$
consisting of
$18 + 6 \lceil3(d-1)/s\rceil
\leq 36$
involutions (elementary matrices in $\EL_3(R)$).
\end{corollary}
\begin{proof}
Using the definition of $\EL_3(R)$ we can see that
if the ring $R$ is generated by $\alpha_k$ then
the group $\EL_3(R)$ is generated by the set
$$
\bar S =
\left\{\Id + e_{i,j} \mid i \not=j \right\}
\cup
\left\{\Id + \alpha_k e_{i,j} \mid i \not=j %,\,\, k=1,\dots,3
 \right\}.
$$
The set $\bar S$ consists of involutions because the ring $R$ has
characteristic $2$.
\end{proof}
\iffalse
\begin{remark}
Actually the group $\Delta$ can be generated by only $10$ involutions.
\end{remark}
\fi

%\pagebreak

The main result in this section is:
\begin{theorem}
\label{KCdelta}
The \KaC\ of the group $\Delta$ with respect to the set
$\bar S$ is
$$
\KC(\Delta; \bar S) \geq \frac{1}{550}.%??
$$
\end{theorem}

\begin{proof}
The proof uses bounded generation of
$\EL_3(R)$ and the following Theorem~2.2 from~\cite{KSL3k}
(see also Theorem~3.4 in~\cite{YSh}):
\begin{theorem}
\label{relTconst}
Let $R$ be an associative  ring generated by $1, \alpha_1,\dots,\alpha_t$.
Let $F_1$ be set of $4(t+1)$ elementary matrices in $\EL_2(R)$ with
$\pm 1$ and $\pm \alpha_i$ off the
diagonal and $F_2$ be the set of standard basis vectors in ${R}^2$
and their inverses. Then the relative \KaC\ of the pair
$(\EL_2(R)\ltimes R^2,R^2)$ is
$$
\KC(\EL_2(R)\ltimes R^2,R^2; F_1 \cup F_2) \geq
\frac{1}{\sqrt{18}(\sqrt{t} + 3)}.
$$
\end{theorem}
\begin{proof}
This is only a short sketch of proof, for details the reader is
referred to~\cite{KSL3k} or~\cite{YSh} in the case of commutative ring.
With out loss of generality we may assume that the ring $R$ is
the free associative
ring generated by $\alpha_i$. For any unitary representation
$\rho : \EL_2(R)\ltimes R^2 \to U(\mr{H})$ and a unit vector $v \in \mr{H}$
we can construct a measure $\mu_v$ on the dual
$\widehat{R^2}$.
%$$
%\widehat{R^2} \simeq \left(\R/\Z [[ \alpha_i^{-1}]] \right)^2.
%$$

If the vector $v$ is an almost invariant under the set $F_1$
then we have
$$
|\mu_v(B) - \mu_v(gB)| \ll 1,
$$
for any Borel set $B \subset \widehat{R^2}$ and any $g \in F_1$.
This show that $\mu_v$ is almost invariant measure on $\widehat{R^2}$.
The vector $v$ is almost invariant for $F_2$ therefore
$\mu_v(B_i) \ll 1 $ for some specific Borel sets $B_i$.
It is known that there are only few $\EL_2(R)$
invariant measures on $\widehat{R^2}$ and almost all of them give large
measures to the sets $B_i$.
It can be shown that the above inequalities imply that
$\mu(\{0\}) >0$.
By definition $\mu(\{0\})$ is square of the length
of the projection of $v$ onto the subspace of $H$ invariant vectors.
Thus $\mu(\{0\}) >0$ implies the existence of $H$ invariant vectors in
$\mr{H}$.
%{\bf Sketch the proof of this theorem}
\end{proof}

Using $6$ different embeddings of $\EL_2(R)\ltimes R^2$ into
$\EL_3(R)$ we can obtain the following
implication of the above theorem:
\begin{equation}
\label{relKC1}
\begin{array}{r@{}c@{}l}
\KC(\Delta;\bar S) & {} \geq {} & \displaystyle
\frac{1}{2}
\KC(\EL_2(R)\ltimes R^2,R^2; F_1 \cup F_2)
\KC(\Delta;GEM) \geq\\
& \geq & \displaystyle
\frac{1}{6\sqrt{2}(3 +\sqrt{5})}\KC(\Delta;GEM)
%\geq
%\frac{1}{42\sqrt{2}}\KC(\Delta;GEM),
\end{array}
\end{equation}
where $GEM$ is the set of all generalized elementary matrices in
$\EL_3(R)$, i.e., the set of all matrices of the form
$$
\left(
\begin{array}{ccc}
1 & 0 & * \\
0 & 1 & * \\
0 & 0 & 1 \\
\end{array}
\right)
\quad
\mbox{or}
\quad
\left(
\begin{array}{ccc}
1 & * & * \\
0 & 1 & 0 \\
0 & 0 & 1 \\
\end{array}
\right)
$$
up to permuting the rows and the columns.

\begin{lemma}
Any element $g$ from the group $\Delta$ can be written as
a product of $17$ elements from the $GEM$, i.e., we have
$\Delta = GEM^{17}$.
\end{lemma}
\begin{proof}
Let $g \in \EL_3(R)$.
With three additional left multiplications by GEMs we can transform $g$
to a $3\times 3$ block matrix where the last column is trivial,
with an extra $3$ left multiplications by GEMs we can make
the second column trivial. Finally with one GEM we can transform $g$
to a matrix which differs form the identity only in the top-left corner,
i.e., we have:
$$
\left(\begin{array}{ccc}
* & * & * \\ * & * & * \\ * & * & *
\end{array} \right)
\stackrel{3}{\Longrightarrow}
\left(\begin{array}{ccc}
* & * & 0 \\ * & * & 0 \\ * & * & 1
\end{array} \right)
\stackrel{3}{\Longrightarrow}
\left(\begin{array}{ccc}
* & 0 & 0 \\ * & 1 & 0 \\ * & 0 & 1
\end{array} \right)
%\stackrel{1}{\Longrightarrow}
\stackrel{1}{\Longrightarrow}
\left(\begin{array}{ccc}
* & 0 & 0 \\ 0 & 1 & 0 \\ 0 & 0 & 1
\end{array} \right).
$$
The entry in the top-left corner is an element in
$\SL_s(\F_2)^{\times K^{d-1}}$ and thus is a group commutator of
two invertible elements in $R$. As such it can be written as a product of
$10$ generalized elementary matrices in $\EL_3(R)$.
This shows that every matrix in $\SL_{3}(R)$ can be written as a product
of $17$ generalized elementary matrices which are in the set $GEM$.
\end{proof}

\begin{remark}
This proof works only for finite rings $R$, which have the property
that any element in $[R^*,R^*]$ is a commutator of two elements
in the multiplicative group $R^*$.
\end{remark}

This together with proposition~\ref{KCball} implies
\begin{corollary}
We have
$$
\KC(\Delta;GEM) \geq \frac{\sqrt{2}}{17}.
$$
\end{corollary}

The above corollary together with equation~(\ref{relKC1}) imply
$$
\KC(\Delta;\bar S) \geq
\frac{1}{6\sqrt{2}(3 + \sqrt{5}} \KC(\Delta;GEM)
\geq \frac{1}{6 (3 +\sqrt{5}) \times 17 } > \frac{1}{550}, %535 works
$$
which finishes the proof of Theorem~\ref{KCdelta}.
\end{proof}

Using the embedding $\pi_i$ and proposition~\ref{relKC} we obtain:
$$
\KC(\Alt(N);\cup\pi_i(\bar S)) >
\frac{1}{1100} \KC(\Alt(N);\cup \pi_i(\Gamma)).
$$

\section{Representations of $\Alt(N)$}
\label{cube}

In this section we will use the same notation as in Section~\ref{relT}:
Let $N=K^d$ for some odd $K$ and $d=6$.
We will think that $\Alt(N)$ acts on
the points in a $d$ dimensional cube of size $K$.
Let $\Gamma$ denote the group
$(\Z/K\Z)^{\times K^{d-1}}$. The group $\Gamma$
has $d$ embeddings $\pi_i$ in $\Alt(N)$.
The image of each cyclic subgroup under $\pi_i$
shifts the points on some line parallel to the
$i$-th coordinate axis. Let $E_i=\pi_i(\Gamma)$ denote
the images of $\Gamma$ and let $E$ be  the union
$$
E = \bigcup \pi_i(\Gamma).
$$

The main result in this section is the following:

\begin{theorem}
\label{KCGammas}
The \KaC\ of $\Alt(N)$ with respect to the set $E$
satisfies
$$
\KC(\Alt(N);E) \geq \frac{1}{70},%???
$$
provided that $K$ is odd and $K>10^6$.% is large enough.
\end{theorem}
\begin{remark}
The proof of Theorem~\ref{KCGammas} gives that
$$
\liminf \KC(\Alt(N);E) \geq \frac{1}{60}
\quad \mbox{as} \quad  K \to \infty
$$
and the proof of this statement is slightly easier
because we can ignore many terms which tend to $0$ as
$K \to \infty$. If $K \ll 10^6$ this method give very weak
bounds for the \KaC, however we believe that
Theorem~\ref{KCGammas} also holds for small $K$.
\end{remark}
\begin{proof}
Let
$\rho: \Alt(N) \to U(\mr{H})$
be a unitary representation
of the alternating group, and let $v \in \mr{H}$
be $\epsilon$-almost invariant unit vector for the set $E$.
We will show that $\mr{H}$ contains an invariant vector if
$\epsilon$ is small enough. Without loss of generality we may
assume that $\mr{H}$ is generated by $v$ as an $\Alt(N)$ module.

The irreducible representations of the symmetric group $\Sym(N)$
are parameterized by the partitions $\lambda$ of $N$. Almost the
same is true the representations of the alternating group, but
the correspondence in this case is not 1-to-1, see~\cite{J,JK}.
We are going to avoid this problem by inducing
the representation $\rho$ to the symmetric group
and working with the induced representation
$$
\rho^s = \Ind_{\Alt(N)}^{\Sym(N)} \rho
\quad \quad
\rho^s : \Sym(N) \to U(\mr{H}^s).
$$
Without loss of generality we can assume that $\mr{H} \subset \mr{H}^s$,
then $\mr{H}^s$ is also generated, as $\Sym(N)$-module, by the vector $v$.
We can decompose $\mr{H}^s$ as
$$
\mr{H}^s = \bigoplus_\lambda \mr{H}_\lambda,
$$
where $\mr{H}_\lambda$ is the sum of all irreducible components in
$\mr{H}$ which correspond to the partition $\lambda$.
We can group these terms in three parts
and break $\rho^s$ as a sum of three representations
$$
\mr{H}^s = \mr{H}_1 \oplus \mr{H}_2 \oplus \mr{H}_3,
$$
where
$$
\mr{H}_1 =
\bigoplus_{
\renewcommand{\arraystretch}{.5}
\begin{array}{c}
{\scriptscriptstyle \lambda'_1< N - h}\\
{\scriptscriptstyle \lambda_1< N - h}
\end{array}
}
\mr{H}_\lambda
\quad
\mr{H}_2 = \bigoplus_{\lambda_1\geq N - h} \mr{H}_\lambda
\quad
\mr{H}_3 = \bigoplus_{\lambda'_1\geq N - h} \mr{H}_\lambda
$$
where $\lambda'$ denotes the dual partition and
$h=\frac{1}{2}K^{3/2}$. %, for $\beta=\frac{7}{5}$.
The action of $\Alt(N)$ on $\mr{H}_\lambda$ and $\mr{H}_{\lambda'}$
is the same for any partition $\lambda$, thus
$\mr{H}_2$ and $\mr{H}_3$ are isomorphic as
representations of the alternating group. Without
loss of generality we may assume that the unit vector $v=v_1+v_2$ has
components only in $\mr{H}_1$ and $\mr{H}_2$.

In this section we will use several probability arguments,
which are based on the following lemma:
\begin{lemma}
\label{proba}
Suppose there we have $lk$ urns grouped in $l$ boxes of size $k$.
If we put $p$ balls in these urns such that different balls go to
different urns, then the probability of having
at least  $q$ balls in the first box is less than
$$
P(l,k,p,q) =
\left(\begin{array}{@{}c@{}} p \\ q \end{array} \right)
\left(\frac{k}{lk-p}\right)^q
\leq
\left(\frac{ep}{ql}\right)^q
\exp\left( \frac{qp}{kl-p}\right).
%\left(1-\frac{p}{kl}\right)^{-q}.
$$
If we drop the restriction that different balls go to different urns then
we can omit the exponential factor from the above estimate.
\end{lemma}
\begin{proof}
Since different balls have to go to different urns, the relative probability
that some ball ends in the first box is not exactly $1/l$, but it is
clear that this relative  probability is
less than $\frac{k}{lk-p}$. This allows us to bound
the number of balls in the first box using the binomial distribution
with $\lambda=\frac{k}{lk-p}$.
Let
$$
B_{p,\lambda}(x) = (1-\lambda + \lambda x)^p = \sum b_ix^i
$$ denote the generating function
of the binomial distribution with parameters $\lambda$ and $p$.
It is clear that the probability of having at least
$q$ balls in the first box is less than
$$
P=\sum_{i\geq q} b_i.
$$
We have that
$$
\frac{1}{q!} \frac{d^qB_{p,\lambda}(x)}{dx^q} =
\sum_{i \geq q} \left(\begin{array}{@{}c@{}} p \\ q \end{array} \right)
b_i x^{i-q} \geq \sum_{i \geq q} b_i x^{i-q},
$$
i.e.,
$P \leq \frac{1}{q!} B_{p,\lambda}(1)^{(q)}$.
The explicit formula for $B_{p,\lambda}(x)$ gives us
$$
P \leq
\left.
\frac{1}{q!} \frac{d^q(1-\lambda + \lambda x)^p}{dx^q}
\right|_{x=1}
= \left(\begin{array}{@{}c@{}} p \\ q \end{array} \right) \lambda^q,
$$
which implies that
$$
P(l,k,p,q) =
\left(\begin{array}{@{}c@{}} p \\ q \end{array} \right)
\left(\frac{k}{lk-p}\right)^q.
$$
The second inequality follows from the estimates
$$
\left(\begin{array}{@{}c@{}} p \\ q \end{array} \right)
\leq \frac{p^q}{q!} \leq \frac{p^q}{(q/e)^q} = \left(\frac{ep}{q}\right)^q
$$
and
$$
\frac{lk}{lk-p} = 1 + \frac{p}{lk-p} < \exp\left(\frac{p}{lk-p}\right).
$$
\end{proof}

Next we prove the following theorem, which together
with estimates for the values of the characters of the
symmetric group from~\cite{Ro}, implies that the component $v_1$ of $v$
is short:
\begin{theorem}
\label{congclass}
The ball $E^{47}$ of radius $47$ generated by the set $E$ contains
almost all cycles in $\Alt(N)$ of length
approximately $K^{d-1}/3\ln K$, %where $\alpha=\frac{9}{5}$,
i.e.,
$$
\left| C_L \setminus E^{47} \right| \leq \frac{1}{100} \left| C_L \right|,% ???
$$
where $C_L$ is the conjugancy class in $\Alt(N)$ of all cycles of
length $L$, which is the largest integer of the from
$1 + a(K-1)$ less than $K^{d-1}/3\ln K$.%
%\footnote{The same argument works for any $L \leq K^2/50$.}
\end{theorem}
\begin{proof}
Let $M$ denotes the set of points in the $d$ dimensional cube with first
coordinate equal to $1$.

\begin{lemma}
\label{tosquare}
Let us chose $L$ distinct points in the $d$ dimensional cube of size $K$.
With high probability (more than
$.99$)
we can move
these points to the set $M$ using only $2$ elements from $E$,
provided that $K\geq 10^6$.
\end{lemma}
\begin{proof}
The probability of having at least $K/\ln K$ points, which
differ only in the first or second coordinate is less than
$$
{\renewcommand{\arraystretch}{2.5}
\begin{array}{r@{}c@{}l}
P_1 & {} = {} &
\displaystyle
K^{d-2} \times
P\left(K^{d-2},K^2,\frac{K^{d-1}}{3\ln K} ,\frac{K}{\ln K}\right)
\leq \\
&\leq & \displaystyle
K^{d-2}
\left(
\frac{e\frac{K^{d-1}}{3 \ln K}}{K^{d-2} \times \frac{K}{\ln K}}
\right)^{\frac{K}{\ln K}}
\exp\left(\frac{\frac{K^{d-1}}{3\ln K} \times \frac{K}{\ln K}
}{K^{d} - \frac{K^{d-1}}{3\ln K}}\right)
\leq \\
& \leq & \displaystyle
K^{d-2} \left(\frac{e}{3}\right)^{\frac{K}{\ln K}}
\exp\left( \frac{1}{10}\right)
\leq \exp(-70),
%\\
%& \leq & \displaystyle
%K \exp\left(\frac{\left(3.9 + 20(\alpha-2)\right) K}{20}\right) e
%\leq \\
%& \leq & \displaystyle
%K\exp\left(-\frac{K}{200}\right)e \leq \exp\left(-e^{10}\right).
\end{array}
}
$$
where the finial inequality holds only if $K$ is large enough.
Also the probability of having at least  $\ln K$ points which
differ only in the second coordinate is less than
$$
{\renewcommand{\arraystretch}{2.5}
\begin{array}{r@{}c@{}l}
P_2 & {} = {} &
\displaystyle
K^{d-1} \times P \left(K^{d-1},K,\frac{K^{d-1}}{3\ln K}, \ln K\right)
\leq \\
& \leq & \displaystyle
K^{d-1}
\left(
\frac{e\frac{K^{d-1}}{3\ln K}}{K^{d-1} \times \ln K }
\right)^{\ln K}
\exp\left(
\frac{\frac{K^{d-1}}{3\ln K} \times  \ln K}{K^{d} -\frac{K^{d-1}}{3\ln K}}
\right)
\leq \\
& \leq & \displaystyle
K^{d-1}
\left(
\frac{e}{3(\ln K)^2}\right)^{\ln K}
\exp\left( \frac{1}{2K}
\right)
\leq
K^{-0.35}
%\frac{e^{20}}{K^2}\left(1-K^{\alpha-2}\right)^{-20}
%\leq
%e^{-15} - \exp(-e^{10}).
\end{array}
}
$$
The above inequalities imply that almost surely (with probability more than $1-P_1-P_2> .99$)
there are no more than $\ln K$ points in
each line parallel to the second coordinate axis and no more
than $K/\ln K$ points in each square parallel to the first and
second coordinate axis.

We claim that if we have a set $b$ of $L$ points in
such good position then
there exists $g \in E_2$ such that there are no two points in
$g b$ which differ only in the first coordinate.

Such an element $g$ can be constructed as follows:
Let us enumerate the lines parallel to the second coordinate
axis and let $b_i$ be the subset of $b$ consisting of all
points on the first $i$ lines. We will prove by induction that
there exists $g_i$ such that there are no two points in
$g_ib_i$ which differ only by the first coordinate. The base case
is trivial and $i=K^5$ will prove the claim. The induction step
is the following: The action of $g_{i+1}$ on all the lines but the $i$th
will be the same as the action of $g_i$. There are $K$ possibilities
for the action on the $i$th line. The number of bad choices (the ones where
there are two points in $g_{i+1}b_{i+1}$ which differ only in the
first coordinate)
is at most
$(s_i-t_i) \times t_i$, where $t_i$ is the number of points from $b$
on the $i+1$'st line and $s_i$ is the number of points from $b$ in the square
parallel to the first and second coordinate axis containing the
$i+1$ line. By assumption the points in $B$ are in a
good position, therefore $t_i \leq \ln K$ and $s_i \leq K/\ln K$, thus
$$
(s_i-t_i) \times t_i \leq (s_i-1) t_i \leq (\ln K-1) \times K/\ln K < K.
$$
This shows there is some good choice for the action of $g_{i+1}$ on
the $i+1$'st line. Thus we can define the action of $g_{i+1}$ on the $i+1$'st
line such that no two points in $g_{i+1}b_{i+1}$ differ only in their
first coordinate, which proves the induction step and the claim.

If no two points from $gb$ differ only in the first coordinate,
we can find
$h \in E_1$ such that $h g b$ is a subset of $M$.
%({\bf explain!!!})
\end{proof}

\begin{lemma}
\label{symonsquare}
For any permutation $\sigma \in \Sym(K^{d-1})$ acting on the points
in $M$ there exist an element
$t \in E^{4d-5}$ such that the restriction of $t$ on this square
is the same as $\sigma$.
\end{lemma}
\begin{proof}
Notice that inside $E_1E_iE_1$ there is a subgroup which preserves $M$
and acts on it as $\Sym(K)^{\times K^{d-2}}$. This subgroup
preserves all coordinates but the $i$th one.
Using the butterfly lemma which gives us that any permutation $g$ in
$\Sym(|A|\times |B|)$ can be written as $g = abc$ where
$a,c \in \Sym(A)^{\times |B|}$ and $b \in \Sym(B)^{\times |A|}$,
it is easy to show that there is an element in
$$
E_1E_2E_1E_3E_1 \dots E_{d-1}E_1 E_d E_1 E_{d-1} \dots E_1E_2E_1 \subset E^{4d-5}
$$
which preserves $M$ and acts on it as $\sigma.$
\end{proof}

\begin{lemma}
\label{cycle}
For any integer $a$ less than $\frac{K^{d-1}-1}{K-1}$
there is an element in
$c_0 \in E_2E_3E_4E_5E_6$
which is a cycle of length $1 + a(K-1)$ in the face $M$.
\end{lemma}
\begin{proof}
Chose $a$ lines in $M$ such that each line is parallel to
some coordinate axis and their union is a tree. Then the products
of the shifts on these lines in any order is a cycle on their
union which contains exactly $1 + a(K-1)$ points.
By definition each shift is in some $E_i$, which proves that there is
an element in $E_2\dots E_6$ which acts as a cycle of
length $1 + a(K-1)$.
%Take $c_2\in \pi_2(\Gamma)$ which acts as long cycle in exactly $a$  rows form $M$, and
%$c_3\in \pi_3(\Gamma)$ which acts as long cycle in exactly one column form $M$.
%It is easy to see that $c=c_2c_3$ is a cycle with support the
%union of the supports of $c_2$ and $c_3$.
%{\bf Explain}
\end{proof}

Now we can finish the proof of Theorem~\ref{congclass}. Let $c$ be
any cycle of length $L$. By lemma~\ref{tosquare} almost surely
the support of $c$ can be moved to $M$ by some element in $E^2$.
Lemma~\ref{symonsquare} gives us the $c$ can be conjugated to
$c_0$ by an element in $E^2\times E^{19}= E^{21}$, i.e., that $c \in E^{47}$.
\end{proof}

Having that $E^{47}$ contains almost a whole conjugancy class we can use
the Roichman~\cite{Ro} character estimates: Let $\lambda$ be a partition
of $N$ and  $C_L$ be the conjugancy class in $\Sym(N)$ consisting of
all $L$-cycles. Then the character $\chi_\lambda$ of the
irreducible representation corresponding to the partition $\lambda$
satisfies the inequality:
\begin{equation}
\label{Roes}
|\chi_\lambda(C_L) | \leq \chi_\lambda(id)
\max\left\{
\frac{\lambda_1}{N}, \frac{\lambda'_1}{N}, \frac{3}{4}
\right\}^{\frac{L-5}{4}}.
\end{equation}

Let
$$
C = \frac{1}{|C_L|}\sum_{g\in C_L} \rho^s(g).
$$
By definition we have
$
\displaystyle
Cv_\lambda =
\frac{\chi_\lambda(C_L)}{\chi_\lambda(id)} v_\lambda
$ for any $v_\lambda \in \mr{H}_\lambda$.
By equation~(\ref{Roes}) we have
$$
{\renewcommand{\arraystretch}{2.5}
\begin{array}{r@{}c@{}l}
|| C v_1 || & {}\leq {} & \displaystyle
||v_1||
\max_{
\left\{
\lambda \big|
%\stackrel{\lambda_1< N - h}{\scriptscriptstyle \lambda'_1< N - h}
{\renewcommand{\arraystretch}{.5}
\begin{array}{c}
{\scriptscriptstyle\lambda_1< N - h}
\\
{\scriptscriptstyle \lambda'_1< N - h}
\end{array}
}
\right\}}
\max\left\{
\frac{\lambda_1}{N}, \frac{\lambda'_1}{N}, \frac{3}{4}
\right\}^{\frac{L-5}{4}}
\leq
\\
& \leq & \displaystyle
||v_1|| \left(1 - \frac{h}{N} \right)^{\frac{L-5}{4}}
\leq
||v_1|| \exp\left(- \frac{h(L-5)}{4N}\right)
\leq
\\
& \leq & \displaystyle
||v_1|| \exp\left(- \frac{\frac{1}{2}K^{3/2}(\frac{K^{d-1}}{3\ln K}-K-4)}
{4K^d}\right)
\leq
\\
& \leq & \displaystyle
||v_1|| \exp\left(- \frac{K^{1/2}}{24 \ln K}
%\left(1-\frac{3(K+4)\ln K}{K^{d-1}}\right)\right)
\left(1-3(K+4)K^{1-d}\ln K\right)\right)
\leq
\\
& \leq & \displaystyle
%||v_1|| \exp\left(- \frac{K^{1/2}}{200}\right)
%\leq
%||v_1|| \exp\left(- \frac{K^{\frac{1}{5}}}{5}\right)
%\leq
%\\
%& \leq & \displaystyle
%||v_1|| \exp\left(- \frac{e^4}{5}\right)
%\leq
e^{-3} ||v_1||.
\end{array}
}
$$

On the other hand we have
$$
{\renewcommand{\arraystretch}{2.5}
\begin{array}{r@{}c@{}l}
||Cv_1 - v_1 ||
& {}\leq {}& \displaystyle
||Cv - v || \leq
\frac{1}{|C_L|}\sum_{g\in C_L} ||\rho^s(g)v -v ||
\leq
\\
& \leq & \displaystyle
\frac{1}{|C_L|}\sum_{g\in C_L\cap E^{47} } ||\rho^s(g)v -v ||  +
\frac{1}{|C_L|}\sum_{g\in C_L \setminus E^{47} } ||\rho^s(g)v -v ||
\leq
\\
& \leq & \displaystyle
 \frac{1}{|C_L|}\sum_{g\in C_L\cap E^{47} } 47 \epsilon +
\frac{1}{|C_L|}\sum_{g\in C_L \setminus E^{47} } 2
\leq
\\
& \leq & \displaystyle
47 \epsilon + 2 \frac{|C_L \setminus E^{47}|}{|C_L|} \leq
47 \epsilon + 0.02,
\end{array}
}
$$
because the elements in $E^{47}$ move the vector $v$ by
at most $47 \epsilon$ and the other elements move $v$ by less than $2$.
Combining the above two equations gives us
\begin{equation}
\label{h2}
{\renewcommand{\arraystretch}{1.5}
\begin{array}{r@{}c@{}l}
||v_1|| & {} \leq {} &
\displaystyle
||Cv_1 - v_1|| + ||Cv_1|| \leq 47 \epsilon + 0.02 + e^{-3} ||v_1||
\leq
\\
& \leq &
\displaystyle
47 \epsilon + 0.02 + e^{-3} < 47 \epsilon + 0.07.
\end{array}
}
\end{equation}
\iffalse
which implies that
\begin{equation}
\label{h2}
||v_1|| \leq
\frac{47 \epsilon + 8.10^{-5}}{1-e^{-50}}\leq 48 \epsilon + 10^{-4}.
\end{equation}
\fi

\medskip
In order to show that $v_2$ is close to an invariant vector we
use entirely different argument.

Let $\mr{B}$ be the set of all ordered $h$-tuples of points in the
six dimensional cube
with side  $K$. Let $\mr{L}$ is the Hilbert space with basis $\mr{B}$
with the natural action of $\Sym(N)$ on it.

First we will prove the next lemma which will allow us to
work with the vector space $\mr{L}$.
\begin{lemma}
There is a injective homomorphism of $\Sym(N)$-modules
$$
i : \mr{H}_2 \to \mr{L}.
$$
\end{lemma}
\begin{proof}
Let decompose $\mr{H}_2$ into sum of irreducible representaions
$$
\mr{H}_2 \simeq \bigoplus_{\lambda_1 \geq N -h} \mr{H}_\lambda
\simeq \bigoplus_{\lambda_1 \geq N -h} m_\lambda V_\lambda,
$$
where $V^\lambda$ irreducible representing of $\Sym(N)$ corresponding
to $\lambda$ and $m_\lambda$ is its multiplicity in $\mr{H}$.
By assumption the space $\mr{H}$ is generated by
the vector $v$, therefore the multiplicities $m_\lambda$ are $0$ or $1$.
The representation theory of $\Sym(N)$ gives us that $V_\lambda$ can
be embedded in $\mr{L}$, which completes the proof.
%Use the description of the irreducible representing $V^\lambda$ to show that
%$V^\lambda$ can be embed in $\mr{L}$ if $\lambda_1 \geq N-h$.
%Also $\mr{H}$ is generated by $v$ therefore the multiplicities of
%$V^\lambda$ in $\mr{H}_2$ are no more than $1$.
\end{proof}

We will prove that the projection $v_{0}$ of $v_2$
onto the space of invariant vectors in $\mr{L}$
is close to $v_2$ by show that the random
walk generated by
$E$ mixes rapidly (in a fixed number of steps)
on the basis $\mr{B}$ of $\mr{L}$.
We will reduce the problem to a random walk on the $h$-tuples of
points in the square of size $K^{d/2}$.%
\footnote{Here we use the assumption that $d$ is even.}
\iffalse
\begin{lemma}
There exists labeling of the points such that
the the two images $E_i = \hat \pi_i(\hat \Gamma)$ of group
$\hat{\Gamma} \simeq (\Z/K^3\Z)^{\times K^{3}}$
lie in $E_1E_2E_3$ and $E_4E_5E_6$ respectively.
\end{lemma}
\begin{proof}
Explain how to obtains powers of some long cycle in $K^3$ as products of
exactly $3$ elements from $E$.
\end{proof}
\fi

Define the operators $U_i$, for $i=1,\dots, 6$ on the
space $\mr{L}$ with basis
$\mr{B}$ by
$$
U_i b = \frac{1}{|E_i|}\sum_{g\in E_i} g b.
$$
Notice that $U_i. U_i =U_i$, because
$E_i$ is a subgroup of $\Alt(N)$.
Let $Q_1$ and $Q_2$ denote the operators
$Q_1 = U_1U_2U_3$ and $Q_2=U_4U_5U_6$.

The next theorem is a quantitative version of the
observation that the random walk on
$\mr{B}$
generated by
$E = \cup E_i$
mixes in few steps, % independent on $K$,
provided that  $h < K^{3/2}$.
\begin{theorem}
\label{ave}
The entries of the matrix $\{a_{b,b'}\}$ of the operator
$Q_2Q_1Q_2Q_1$ defined
by
$$
Q_2Q_1Q_2Q_1 b=
\sum_{b'\in \mr{B}} a_{b,b'} b',
$$
satisfy the inequality
$$
a_{b,b'} \geq \frac{1}{|\mr{B}|}\left( 1- \frac{h^2}{K^3} \right).
$$
This inequalities imply that the operator norm of
$Q_2Q_1Q_2Q_1 - P_0$ is less than $1- \frac{h^2}{K^3}$,
where $P_0$ is the projection onto the space of $\Sym(N)$ invariant
vectors in $\mr{L}$.
\end{theorem}
\begin{proof}
Let $\mr{B}_1$ ($\mr{B}_2$) denote
the sets of all $h$-tuples from $\mr{B}$
such that there are no two points with the same first (last)
three coordinates.
\begin{lemma}
The entries of the matrices $\{q_{1;b,b'}\}$ and $\{q_{2;b,b'}\}$
of the operators $Q_1$ and $Q_2$
satisfy
$$
\sum_{b' \in \mr{B}_1} q_{1;b,b'} \geq 1-\frac{h^2}{2K^3}
\quad
\mbox{and}
\quad
\sum_{b' \in \mr{B}_2} q_{2;b',b} \geq 1-\frac{h^2}{2K^3}.
$$
\end{lemma}
\begin{proof}
Let $p_{i,j,k}$ be the number of points form $b$ which
have $4$th, $5$th and $6$th coordinate equal to $i$, $j$, and $k$
respectively. Lets choose some ordering of the triples
$\{i,j,k\}$ and let $p_s$ be the numbers $p_{i,j,k}$ in this order.

Then the number of elements in $E_1E_2E_3$ which sends $b$ to
an element in $\mr{B}_1$ is at least %{\bf explain}
$$
\begin{array}{r@{}c@{}l}
\displaystyle
Q_{1;b}
&{}\geq {} &
\displaystyle
\prod_s \left( K^{3K^2} - (p_s^2/2 + p_s \sum_{t<s} p_t)
\times K^{3K^2 -3}  \right)
=
\\
&= &
\displaystyle
K^{3K^5} \prod_s \left( 1 - \frac{p_s^2/2+ \sum_{t<s} p_s}{K^3} \right)
\geq
\\
& \geq &
\displaystyle
K^{3K^5} \left( 1 - \frac{\sum_s p_s^2/2 + \sum_s p_s \sum_{t<s} p_t }{K^3}\right)=
\\
& = &
\displaystyle
K^{3K^5}\left( 1 - \frac{\sum_s p_s^2/2 +\sum_{t<s} p_s  p_t }{K^3} \right)
=
\\
& = &
\displaystyle
K^{3K^5}\left( 1 - \frac{\sum_{s} p_s^2 }{2K^3} \right)
=
K^{3K^5}\left( 1 - \frac{h^2}{2K^3} \right)=
\\
& = &
\displaystyle
|E_1| \times |E_2| \times |E_3| \left( 1 - \frac{h^2}{2K^3} \right),
\end{array}
$$
because there are $K^{3K^2}$ possibilities for the action of
$E_1E_2E_3$ on the cube corresponding to $s$, and the number of bad ones
(such that $g b\not \in \mr{B}_1$, because some point from this cube
have the same first $3$ coordinates as one of the previous points)
is at most
$$
\sum_{i=1}^{p_s} (i-1 +\sum_{t<s} p_t) \times K^{3K^2 -3} \leq
\left( \frac{p_s^2}{2} + p_s \sum_{t<s} p_t \right) \times K^{3K^2 -3}.
$$

Using the definition of $Q_1$ and $q_{1;b;b'}$ we can see that
$$
\sum_{b' \in \mr{B}_1} q_{1;b,b'}
=
\frac{Q_{1;b}}{|E_1| \times |E_2| \times |E_3|}
\geq
1-\frac{h^2}{2K^3}.
$$
\end{proof}

\begin{lemma}
If $b \in \mr{B}_1$ and
$b' \in \mr{B}_2$
then the number
$$
P(b,b') =
\frac{
\left|\left\{
g_i\in E_i, \mid
g_1g_2g_3g_4g_5 g_6 b = b'
\right\} \right|
}{
|E_1| \times |E_2| \times |E_3|
\times
|E_4| \times |E_5| \times |E_6|
}
=
K^{-6 h}.
$$
does not depend on $b_1$ and $b_2$ .
%and $P(b_1,b_2) \geq \frac{\gamma}{|\overline{\mr{B}}|}$.
\end{lemma}
\begin{proof}
Let $g_i\in E_i$  satisfy
$g_1 \dots g_6 b=b'$.
If $g_1 \dots g_6$ sends the tuple
$(a_1,\dots,a_6)$ to the tuple $(b_1,\dots,b_6)$, then
$g_6$ acts on the the line $(a_1,\dots,a_5, *)$ as shift by $b_6-a_6$;
$g_5$ acts on the the line $(a_1,\dots,a_4, *,b_6)$ as shift by $b_5-a_5$
and so on.
Thus $g_1 \dots g_6 b=b'$ determines the action of each $g_i$ on
exactly $h$ lines (the conditions $b \in \mr{B}_1$ and $b' \in \mr{B}_2$
imply that all these lines are different).
This shows that
$$
P(b,b') =
\left(
\frac{K^{K^5- h}}{K^{K^5}} \right)^6
\iffalse
\frac{
K^{K^5- h} \times
K^{K^5- h} \times
K^{K^5- h} \times
K^{K^5- h} \times
K^{K^5- h} \times
K^{K^5- h}
}{
K^{K^5} \times
K^{K^5} \times
K^{K^5} \times
K^{K^5} \times
K^{K^5} \times
K^{K^5}
}
\fi
= K^{-6 h}.
$$
\end{proof}

Now we can finish the proof of Theorem~\ref{ave}.
We can write $Q_2Q_1Q_2Q_1 b$ as $Q_2(Q_1Q_2)(Q_1 b)$,
therefore
%Finally the Theorem~\ref{ave} follows easily from the above inequalities and the
%observation
$$
\begin{array}{r@{}c@{}l}
a_{b,b'}
&{}={}&
\displaystyle
\sum_{c,c' \in \mr{B}}
q_{1;b,c}P(c,c')q_{2;c',b'}
\geq
\sum_{c\in \mr{B}_1\,c' \in \mr{B}_2}
q_{1;b,c}P(c,c')q_{2;c',b'}
=
\\
& = &
\displaystyle
K^{-6h}\sum_{c\in \mr{B}_1\,c' \in \mr{B}_2}
q_{1;b,c}q_{2;c',b'}
=
K^{-6h}
\sum_{c\in \mr{B}_1} q_{1;b,c}
\sum_{c'\in \mr{B}_1} q_{1;b,c}
\geq
\\
& \geq &
\displaystyle
K^{-6h} \left(1 -\frac{h^2}{2K^3}\right)^2
\end{array}
$$
The size of the basis $\mr{B}$ is
$$
\begin{array}{r@{}c@{}l}
|\mr{B}|
& {} = {} &
\displaystyle
\prod_{i=1}^h \left(K^6 - i+1\right) =
K^{6h}\prod_{i=1}^h \left(1 - \frac{i+1}{K^6}\right) \geq
\\
& \geq &
\displaystyle
K^{6h}\left( 1 - \sum_{i=1}^h \frac{i+1}{K^6} \right)
\geq
K^{6h}\left( 1 - \frac{h^2}{2K^6} \right).
\end{array}
$$
Theorem~\ref{ave}, follows from the above estimates and the inequality
$$
\left(1 -\frac{h^2}{2K^3}\right)^2 \left( 1 - \frac{h^2}{2K^6} \right)
\geq  \left(1 -\frac{h^2}{K^3}\right).
$$
\end{proof}

Let $v_{0}$ be the projection of $v_1$ on the space of
invariant vectors in $\mr{L}$. Using Theorem~\ref{ave}
we can see
$$
||Q_2Q_1Q_2Q_1v_1 - v_{0} || =
||Q_2Q_1Q_2Q_1(v_1 - v_{0}) ||
\leq \frac{h^2}{K^3} || v_1 - v_{0}||,
$$
but we also have
$$
%||Q_2Q_1Q_2Q_1v_1 - v_{0} || =
||Q_2Q_1Q_2Q_1v_1 - v_1 ||
\leq ||Q_2Q_1Q_2Q_1v - v || \leq 12\epsilon,
$$
which implies that
$$
|| v_1 - v_{0} || \leq 12\epsilon +
\frac{h^2}{K^3} || v_1 - v_{0}||.
$$
Substituting $h= \frac{1}{2}K^{3/2}$ in the above inequality yields
\begin{equation}
\label{h1}
|| v_1 - v_{0} || \leq \frac{1}{1-\frac{1}{4}} 12\epsilon =16 \epsilon.
\end{equation}

This finishes the proof of Theorem~\ref{KCGammas}, because
equations~(\ref{h2}) and~(\ref{h1}) imply that
$$
||v - v_{0}|| \leq ||v_2|| + || v_1 - v_{0} || \leq
47 \epsilon + 0.07 + 16 \epsilon = 63 \epsilon + 0.07.
$$
The last expression is less then $1$ if $\epsilon=1/70$, therefore
$v_{0}$ is not zero and it is an invariant vector in
$\mr{H} \subset \mr{H}^s$
\end{proof}

As mentioned in the introduction Theorem~\ref{altexp6}  follows
immediately from Proposition~\ref{relKC} and
Theorems~\ref{relTconst},~\ref{KCGammas}.

\section{Proof of Theorems~\ref{main} and~\ref{mainsym}}
\label{proofmain}

Theorems~\ref{main} and~\ref{mainsym} follow easily from
Theorem~\ref{altexp6}

\begin{proof} {\bf of Theorem~\ref{main}}
By Theorem~\ref{altexp6} the alternating groups $\Alt(n_s)$
are expanders with respect
to some generating set $F_{n_s}$ for $n_s=\left(2^{3s}-1\right)^6$.
The sequence $\left\{n_s\right\}_{s}$
grows exponentially. Thus, for
any sufficiently large $n>10^{38}$ there exists $s$
such that
$$
1< \frac{n}{n_s} < \max_s \frac{(2^{3s+3}-1)^6}{(2^{3s}-1)^6}
\leq \max_s\left(8+\frac{7}{2^{3s}-1}\right)^{6} < 3 \times 10^5.
$$

Using the butterfly lemma it can be shown that
the group $\Alt(n)$ can be written as a product of
several  copies of $\Alt(n_s)$
embedded in $\Alt(n)$. The number of copies $P$ is at most
$$
P \leq 3 \lceil n/n_s \rceil +3 < 10^6
$$

Let $\pi_i$ denote the $P$ embeddings of $\Alt(n_s)$ in
$\Alt(n)$
Taking the union of $\pi_i(F_{n_s})$ one
obtains a generating set $F_n$ of $\Alt(n)$.
By Propositions~\ref{KCball} and~\ref{KCrel} we have
$$
\begin{array}{r@{}c@{}l}
\displaystyle
\KC\left(\Alt(n); F_n\right) & {}\geq {}&
\displaystyle
\frac{1}{2} \KC\left(\Alt(n); \cup \pi_i(\Alt({n_s}))\right)\KC\left(\Alt(n_s); F_{n_s}\right)
\geq
\\
& \geq &
\displaystyle
\frac{\sqrt{2}}{P} \KC\left(\Alt(n_s); F_{n_s}\right)
\geq 10^{-12},
\end{array}
$$
because every element in $\Alt(n)$ can be written as a product of no more than
$P$ elements from $\cup \pi_i(\Alt({n_s})$.
By the construction it is clear that the size of $F_n$ is
bounded above by $10^9$.

If $n$ is small then the above argument does not work, because
we can not find $n_s<n$ with the desired properties. But there are only
finitely many $n$ less than $10^{38}$.
It is clear for each such $n$ that there exist
$\epsilon_n>0$ and a generating set $F_n$ with less than $10^9$ elements
such that
$\KC(\Alt(n);F_n) \geq \epsilon_n$.

Thus for each $n$ we have constructed a generating set $F_n$ of the alternating
group $\Alt(n)$ of size at most $L=10^9$ such that
$\KC(\Alt(n);F_n) \geq \epsilon$, where
$\epsilon = \min\{10^{-12}, \inf_n \epsilon_n \} >0$, which completes the
proof of Theorem~\ref{main}.
\end{proof}

\begin{proof} {\bf of Theorem~\ref{mainsym}}
The alternating group $\Alt(n)$ is a subgroup of index $2$ inside
$\Sym(n)$. Let $t \in \Sym(n)\setminus \Alt(n)$ be an odd permutation.
If $F_n$ is an expanding generating set of $\Alt(n)$ then
$\tilde F_n =F_n \cup \{t\}$ is an expanding generating set of $\Sym(n)$ and the
Kazhdan constants are almost the same.

%Let $\rho:\Sym(n) \to U(\mr{H})$ is a unitary representation of
%the symmetric group and let $v\in \mr{H}$ be an $\epsilon$ almost
%invariant vector for $\tilde F_n$.
Using the trivial inequality
$$
\KC(\Sym(n),\Alt(n);F_n \cup \{t\}) \geq
\KC(\Alt(n);F_n)
$$
and proposition~\ref{KCrel} we can see that
$$
\KC(\Sym(n);F_n \cup \{t\}) \geq
\frac{1}{2}
\KC(\Alt(n);F_n)\KC(\Sym(n);\Alt(n) \cup \{t\}).
$$
However it is clear that any element in $\Sym(n)$ can be written
as a product of two elements from $\Alt(n) \cup \{t\}$, thus
$$
\KC(\Sym(n);\Alt(n) \cup \{t\}) \geq \frac{\sqrt{2}}{2}.
$$
These two inequalities give us
$$
\KC(\Sym(n);F_n \cup \{t\}) \geq
\frac{\sqrt{2}}{4}
\KC(\Alt(n);F_n)
\geq \frac{1}{3}
\KC(\Alt(n);F_n),
$$
which finishes the proof of Theorem~\ref{mainsym}.
\end{proof}

\begin{remark}
In the above proofs of Theorems~\ref{main} and~\ref{mainsym},
we have obtained huge generating sets $F_n$ and $\tilde F_n$
with very weak bounds for the \KaC s. Using a more careful
analysis it is possible to decrease the size of $F_n$ and to
improve the bounds for the Kazhdan constants.
We can prove that for all
sufficiently large $n$
%\linebreak $n>10^{40}$
there exist generating sets
$F_n$ and $\tilde F_n$ of the alternating and the symmetric groups
such that
\begin{itemize}
\item
$|F_n| \leq 20$ and $|\tilde F_n| \leq 20$,
\item
both $F_n$ and $\tilde F_n$ consist only of involutions,
\item
$\KC(\Alt(n);F_n) \geq 10^{-7}$ and $\KC(\Sym(n);\tilde F_n) \geq 10^{-7}$.
\end{itemize}
\end{remark}

\begin{remark}
The bounds for the \KaC s in Theorems~\ref{main} and~\ref{mainsym}
are explicit in the case $n>10^{40}$. This restriction comes form the
condition $K \geq 10^6$ in Theorem~\ref{KCGammas}.
We suspect that some modification in the proof will allow us
to weaken requirement.
\end{remark}

\section{Comments}
\label{comments}

In this section we will briefly discuss some variations of the construction.

As in sections~\ref{relKC} and~\ref{altexp6}, we will
denote $N=K^d$, where $d$ is fixed.
Let $H$ be any group which acts transitively on the set of $K$
points, we will also assume that all elements of $H$ act as
even permutations, i.e., we have $H \hookrightarrow \Alt(K)$.%
\footnote{If the last requirement is not satisfied then it is
impossible to make the Cayley graphs of groups $\Delta(H)$ expanders
with respect to any generating set of bounded size.}

We will define the group
$$
\Delta(H) = H^{\times K^{d-1}}.
$$
This group can be embedded in $\Alt(N)$ in $d$ different ways
using $\pi_i$. Let $E= \cup \pi_i(\Delta(H))$ be the union of the
images of these embeddings.

The main result in section~\ref{altexp6}, Theorem~\ref{KCGammas},
says that
$$
\KC\left(\Alt(N);\cup \pi_i(\Delta(\Z/K\Z)\right) \geq \frac{1}{60},
$$
provided that $d=6$ and $K$ is large enough. This result can be
generalized to any transitive group:

\begin{theorem}
\label{KCGammas-gen}
If $H$ is a transitive group acting on $K$ points then
$$
\KC(\Alt(N);\cup \pi_i(\Delta(H)) \geq \frac{1}{60},
$$
provided that $d=6$ and $K$ is large enough.
\end{theorem}
\begin{proof}
The proof is almost the same as the one of Theorem~\ref{KCGammas}.
The only difference is that Lemma~\ref{cycle} does not hold. However
there is a weaker analog of this lemma:

\begin{lemma}
\label{cycle-gen}
For any integer $a$ less than $K^{d-2}$ there is an element $c_0$ in $E_2$,
which acts on $M$ as a permutation with $Ka$ non-fixed points.
\end{lemma}
\begin{proof}
Let $g \in H$ be an element without fixed points ---
such an element exists in any transitive group because the average
number of fixed points is $1$. The $c_0$ is the image of an element in
$\Delta(H)$ which is equal to $g$ in $a$ copies of $H$ and is identity in the other
copies.
\end{proof}

The same argument as in the proof of Theorem~\ref{KCGammas} gives that
$E^{47}$ (even $E^{43}$) contains almost all elements in the conjugancy class $C$, which
contains $c_0$. By construction the permutations in this
conjugancy class contain $L$ non-fixed points, where $L$ is
approximately $K^{d-1}/3\ln K$.

For such conjugancy classes there are character estimates, see~\cite{Ro},
similar to~(\ref{Roes}):
\begin{equation}
\label{Roes-gen}
|\chi_\lambda(C) | \leq \chi_\lambda(id)
\max\left\{
\frac{\lambda_1}{N}, \frac{\lambda'_1}{N}, q
\right\}^{cL},
\end{equation}
where $c$ and $q$ are some universal constants. The precise values
of these constants are not in the literature but using the proof of
Theorem~1 from~\cite{Ro} one can obtain estimates for these constants.
Our computations show that
we can use $q=1-10^{-3}$ and $c=10^{-3}$.

The estimate in~(\ref{h2}) continues to hold if $K$ is sufficiently
large. The bound depends on the values of the constants $c$ and $q$ above.
A careful computation of the constants $c$ and $q$ shows that
we can use $q=1-10^{-3}$ and $c=10^{-3}$, which gives that if $K>10^{10}$
then $K$ is sufficiently large and the estimate~(\ref{h2}) holds.
The rest of the proof is the same as in Theorem~\ref{KCGammas}.
\end{proof}

%\begin{remark}
Actually we have shown that
$$
\liminf_{K\to \infty} \KC(\Alt(N);\cup \pi_i(\Delta(\Z/K\Z)) \geq \frac{1}{C(d)},
$$
if $d=6$ and $C(6) = 60$.
This result is also valid for larger values of $d$ with $C(d)=10d$
and the proof is essentially the same.

\medskip

It is interesting to see if a similar result holds for small values of $d$:
In the case $d=1$ the question does not make sense because
the group $\Alt(N)$ is not generated by the set $E$.
If $d=2$, we think that it is not possible to obtain a bound
for the \KaC\ which is independent on $K$, because the group
$\Delta(H)$ is not large enough.

Using more careful analysis of the random walk generated by $E$ on the set
of $h$ tuples, it is possible to show that there is a uniform bound
for the \KaC\ if $d=4$ or $d=5$. We do not know whether  such result is also
valid in the case $d=3$.
%, for details see~\cite{Krwc}.
%\end{remark}

\bigskip

In order to use the groups $\Delta(H)$ to construct  bounded
degree expanders we need to find a family of groups $H_s$
which acts transitively on $K_s$ points such that there is
a uniform lower bound for the \KaC s:
$$
\KC(H_s^{\times K_s^{d-1}}; S_s)
$$
for some generating sets $S_s$.
In our proof we used the groups $H_s=\SL_{3s}(\F_2)$
acting on $F_2^{3s} \setminus\{0\}$.
A more natural family of groups with these properties is
$\SL_{n}(\F_{q})$ acting on
${\F_{q}}^{n} \setminus \{0\}$ or on the projective space  $P^{n-1}\F_{q}$,
where $n\geq 3$ is fixed and $q$ is a power of a prime number.
In the case $n=4$ and $q=2^s$ we can even find an element in
$\SL_4(\F_{2^s})$ which acts as a long cycle on $P^3\F_{2^s}$.
In this case we will  obtain slightly better \KaC s
with respect to some generating set of $\Alt(N)$ where
$N=\left(2^{3s}+2^{2s} + 2^s + 1\right)^6$.

Using the groups $\SL_{3s}(\F_2)$ has some advantages: it is possible
to generate the product of $2^{s^2}$ copies of the matrix algebra
$\Mat_s(\F_2)$ by just three elements. This allows us to
construct generating sets $\bar S_s$ of
$\Delta(\SL_{3s}(\F_2))^{\times s}$ such
that
$$
\KC(\Delta(\SL_{3s}(\F_2))^{\times s}; \bar S_s) \geq \frac{1}{1000}.
$$
Using these generating sets we can make the Cayley graphs of the product of $s$ copies of
$\Alt(N)$ expanders.

Another advantage is the following:
Let $\tilde \Delta$ be the infinite product of the groups $\Delta(\SL_{3s}(\F_2))$.
Inside $\tilde \Delta$ there is a dense subgroup $\bar \Delta$ generated
by the set $S$, which projects to $S_s$ in every factor.
Using the methods from~\cite{KSL3k},
it can be shown that the the pro-finite completion of the group $\bar \Delta$
is slightly larger than $\tilde \Delta$ and it has property $\tau$.
We can use this group to obtain a dense subgroup $G$ inside
$$
\prod_s \Alt \left((2^{3s}-1)^6\right)
$$
which also has property $\tau$. $G$ is the first example of a dense subgroup in
the product of infinitely many alternating groups which does not
map onto $\Z$, for more details see~\cite{KtauSGG} and for other examples of
dense subgroups in such products see~\cite{Pyb}.

\medskip

Theorem~\ref{main} can be viewed as a major stop towards proving the
conjectured suggested by Alex Lubotzky~\cite{Lpr}:%
\footnote{This conjecture was around for a long time, however it had
never appeared in writing because the affirmative answer was not
known even in many simple cases.}
\begin{conjecture}
\label{FSGE}
Let $G_i$ be the family of all
non-abelian finite simple groups. There exists a generating sets $S_i$
(with uniformly bounded size) such that the Cayley graphs
$\mr{C}(G_i,S_i)$ form a family of $\epsilon$-expanders for
some fixed $\epsilon>0$.
\end{conjecture}

A strong supporting evidence for this conjecture is the following well
known fact (see~\cite{expanderbook,lu,YSh}):
The groups of a fixed Lie type
over different finite fields form an expander family, provided that
the rank of the Lie group is at least $2$. However both the size of the
generating set and the expanding constant depend on the rank.
This shows that almost all non-abelian finite simple groups can be
put into infinitely many families such the groups in each family can
be made expanders.
Another supporting evidence for this conjecture is
that for any non-abelian finite simple
group there exist $4$ generators such that the
diameter of the corresponding Cayley graph is logarithmic in
the size of the group (see~\cite{BLK} and~\cite{KR}).
However in these examples it is known that the Cayley graphs
are not expanders.

As with many similar results,
one expects that the proof of Conjecture~\ref{FSGE}
will use the classification of the finite simple groups.

The first major step towards proving Conjecture~\ref{FSGE} was made
in~\cite{KSL3k} --- there it is shown that the Cayley graphs of
$\SL_n(\F_q)$ for any prime power $q$ and infinitely many $n$ can be
made expanders simultaneously by choosing a suitable generating
sets. This can be generalized to all families of finite simple
groups of Lie type of rank at least $2$.

The results in this paper prove the Lubotzky conjecture in case
of the alternating groups, which was believed that this is the
most difficult case. These results, together with some new ones in the
rank one case, are combined in~\cite{KNFS}, which
almost proves Conjecture~\ref{FSGE}.

\providecommand{\bysame}{\leavevmode\hbox to3em{\hrulefill}\thinspace}
\providecommand{\MR}{\relax\ifhmode\unskip\space\fi}
\providecommand{\MRhref}[2]{%
  \href{http://www.ams.org/mathscinet-getitem?mr=#1}{#2}
}
\providecommand{\href}[2]{#2}

\texttt{\\Martin Kassabov, \\
Cornell University, Ithaca, NY 14853-4201, USA. \\
\emph{e-mail:} kassabov@math.cornell.edu }


\begin{thebibliography}{10}

\bibitem{A}
Mikl{\'o}s Ab{\'e}rt, \emph{Symmetric groups as products of abelian subgroups},
  Bull. London Math. Soc. \textbf{34} (2002), no.~4, 451--456. \MR{MR1897424
  (2002m:20006)}

\bibitem{ALW}
Noga Alon, Alexander Lubotzky, and Avi Wigderson, \emph{Semi-direct product in
  groups and zig-zag product in graphs: connections and applications (extended
  abstract)}, 42nd IEEE Symposium on Foundations of Computer Science (Las
  Vegas, NV, 2001), IEEE Computer Soc., Los Alamitos, CA, 2001, pp.~630--637.
  \MR{MR1948752}

\bibitem{BHKLS}
L.~Babai, G.~Hetyei, W.~M. Kantor, A.~Lubotzky, and {\'A}.~Seress, \emph{On the
  diameter of finite groups}, 31st Annual Symposium on Foundations of Computer
  Science, Vol.\ I, II (St.\ Louis, MO, 1990), IEEE Comput. Soc. Press, Los
  Alamitos, CA, 1990, pp.~857--865. \MR{MR1150735}

\bibitem{BLK}
L.~Babai, W.~M. Kantor, and A.~Lubotzky, \emph{Small-diameter {C}ayley graphs
  for finite simple groups}, European J. Combin. \textbf{10} (1989), no.~6,
  507--522. \MR{MR1022771 (91a:20038)}

\bibitem{BaHarpe}
Roland Bacher and Pierre de~la Harpe, \emph{Exact values of kazhdan constants
  for some finite groups}, Journal of Algebra (1994), no.~163, 495--515.
  \MR{MR1262716 (95b:20018)}

\bibitem{Bur}
M.~Burger, \emph{Kazhdan constants for {${\rm SL}(3,{\Z})$}}, J. Reine Angew.
  Math. \textbf{413} (1991), 36--67. \MR{MR1089795 (92c:22013)}

\bibitem{CK}
David Carter and Gordon Keller, \emph{Bounded elementary generation of {${\rm
  SL}\sb{n}({\cal O})$}}, Amer. J. Math. \textbf{105} (1983), no.~3, 673--687.
  \MR{MR704220 (85f:11083)}

\bibitem{CK1}
\bysame, \emph{Elementary expressions for unimodular matrices}, Comm. Algebra
  \textbf{12} (1984), no.~3-4, 379--389. \MR{MR737253 (86a:11023)}

\bibitem{GG2}
Ofer Gabber and Zvi Galil, \emph{Explicit constructions of linear size
  superconcentrators}, 20th Annual Symposium on Foundations of Computer Science
  (San Juan, Puerto Rico, 1979), IEEE, New York, 1979, pp.~364--370.
  \MR{MR598118 (83g:68097)}

\bibitem{Gilman}
Robert Gilman, \emph{Finite quotients of the automorphism group of a free
  group}, Canad. J. Math. \textbf{29} (1977), no.~3, 541--551. \MR{MR0435226
  (55 \#8186)}

\bibitem{J}
G.~D. James, \emph{The representation theory of the symmetric groups}, Lecture
  Notes in Mathematics, vol. 682, Springer, Berlin, 1978. \MR{MR513828
  (80g:20019)}

\bibitem{J2}
\bysame, \emph{Representations of general linear groups}, London Mathematical
  Society Lecture Note Series, vol.~94, Cambridge University Press, Cambridge,
  1984. \MR{MR776229 (86j:20036)}

\bibitem{JK}
Gordon James and Adalbert Kerber, \emph{The representation theory of the
  symmetric group}, Encyclopedia of Mathematics and its Applications, vol.~16,
  Addison-Wesley Publishing Co., Reading, Mass., 1981. \MR{MR644144
  (83k:20003)}

\bibitem{K}
Martin Kassabov, \emph{{Kazhdan Constants for $\SL_n(\Z)$}},
  \texttt{arXiv:math.GR/0311487}.

\bibitem{KtauSGG}
\bysame, \emph{{Property Tau and Subgroup Growth}}, in preparation.

\bibitem{Ksym}
\bysame, \emph{{Symmetric Groups and Expanders}},
  \texttt{arXiv:math.GR/0503204}.

\bibitem{KSL3k}
\bysame, \emph{{Universal lattices and unbounded rank expanders}},
  \texttt{arXiv:math.GR/0502237}.

\bibitem{KNFS}
Martin Kassabov, Alex Lubotzky, and Nikolay Nikolov, \emph{{Finite simple
  groups and expanders}}, in preparation.

\bibitem{KNtau}
Martin Kassabov and Nikolay Nikolov, \emph{{Universal lattices and Property
  Tau}}, \texttt{arXiv:math.GR/0502112}.

\bibitem{KR}
Martin Kassabov and Tim~R. Riley, \emph{{Diameters of Cayley graphs of
  $\SL_n(\Z/k\Z)$}}, \texttt{arXiv:math.GR/0502221}.

\bibitem{kazhdan}
D.~A. Ka{\v{z}}dan, \emph{On the connection of the dual space of a group with
  the structure of its closed subgroups}, Funkcional. Anal. i Prilo\v zen.
  \textbf{1} (1967), 71--74. \MR{MR0209390 (35 \#288)}

\bibitem{Kl}
Maria Klawe, \emph{Limitations on explicit constructions of expanding graphs},
  SIAM J. Comput. \textbf{13} (1984), no.~1, 156--166. \MR{MR731033
  (85k:68077)}

\bibitem{lubwiess}
A.~Lubotzky and B.~Weiss, \emph{Groups and expanders}, Expanding graphs
  (Princeton, NJ, 1992), DIMACS Ser. Discrete Math. Theoret. Comput. Sci.,
  vol.~10, Amer. Math. Soc., Providence, RI, 1993, pp.~95--109. \MR{MR1235570
  (95b:05097)}

\bibitem{Lpr}
Alexander Lubotzky, \emph{private communication}.

\bibitem{expanderbook}
\bysame, \emph{Discrete groups, expanding graphs and invariant measures},
  Progress in Mathematics, vol. 125, Birkh\"auser Verlag, Basel, 1994.
  \MR{MR1308046 (96g:22018)}

\bibitem{lu}
\bysame, \emph{Cayley graphs: eigenvalues, expanders and random walks}, Surveys
  in combinatorics, 1995 (Stirling), London Math. Soc. Lecture Note Ser., vol.
  218, Cambridge Univ. Press, Cambridge, 1995, pp.~155--189. \MR{MR1358635
  (96k:05081)}

\bibitem{LP}
Alexander Lubotzky and Igor Pak, \emph{The product replacement algorithm and
  {K}azhdan's property ({T})}, J. Amer. Math. Soc. \textbf{14} (2001), no.~2,
  347--363 (electronic). \MR{MR1815215 (2003d:60012)}

\bibitem{LubZuk}
Alexander Lubotzky and A.~{\.Z}uk, \emph{On property $\tau$}.

\bibitem{Mar}
G.~A. Margulis, \emph{Explicit constructions of expanders}, Problemy Pereda\v
  ci Informacii \textbf{9} (1973), no.~4, 71--80. \MR{MR0484767 (58 \#4643)}

\bibitem{MW}
Roy Meshulam and Avi Wigderson, \emph{Expanders in group algebras},
  Combinatorica \textbf{24} (2004), no.~4, 659--680. \MR{MR2096820}

\bibitem{Pyb}
L{\'a}szl{\'o} Pyber, \emph{Groups of intermediate subgroup growth and a
  problem of {G}rothendieck}, Duke Math. J. \textbf{121} (2004), no.~1,
  169--188. \MR{MR2031168 (2004k:20056)}

\bibitem{RVW}
Omer Reingold, Salil Vadhan, and Avi Wigderson, \emph{Entropy waves, the
  zig-zag graph product, and new constant-degree expanders and extractors
  (extended abstract)}, 41st Annual Symposium on Foundations of Computer
  Science (Redondo Beach, CA, 2000), IEEE Comput. Soc. Press, Los Alamitos, CA,
  2000, pp.~3--13. \MR{MR1931799}

\bibitem{Ro}
Yuval Roichman, \emph{Upper bound on the characters of the symmetric groups},
  Invent. Math. \textbf{125} (1996), no.~3, 451--485. \MR{MR1400314
  (97e:20014)}

\bibitem{Ro1}
\bysame, \emph{Expansion properties of {C}ayley graphs of the alternating
  groups}, J. Combin. Theory Ser. A \textbf{79} (1997), no.~2, 281--297.
  \MR{MR1462559 (98g:05070)}

\bibitem{RSW}
Eyal Rozenman, Aner Shalev, and Avi Wigderson, \emph{A new family of cayley
  expanders},  (2004).

\bibitem{YSh}
Yehuda Shalom, \emph{Bounded generation and {K}azhdan's property ({T})}, Inst.
  Hautes \'Etudes Sci. Publ. Math. (1999), no.~90, 145--168 (2001).
  \MR{MR1813225 (2001m:22030)}

\end{thebibliography}
\end{document}